\documentclass[twocolumn]{autart}
\usepackage{mathrsfs}
\usepackage{epsfig}
\usepackage{times}
\usepackage{amssymb}
\usepackage{graphicx}
\usepackage[longnamesfirst,round]{natbib}
\usepackage{natbib}
\graphicspath{{figure/}}
\usepackage[none]{hyphenat}
\usepackage[tbtags]{amsmath}
\usepackage{algorithm,algorithmic}
\usepackage{bm}
\usepackage{url}
\usepackage{subfigure}
\usepackage{enumerate}
\usepackage[monochrome]{color}  

\usepackage{color}
\definecolor{gray}{RGB}{128,128,128}

\newcommand{\blue}[1]{{\color{blue} #1}}

\def\Z{\mathbb{Z}}
\newcommand{\R}{\mathbb{R}}

\DeclareMathOperator{\diag}{diag}
\DeclareMathOperator{\rank}{rank}

\newcommand{\hb}{\mathbf{h}}

\newcommand{\cb}{\mathbf{c}}
\newcommand{\ub}{\mathbf{u}}
\newcommand{\vb}{\mathbf{v}}

\newcommand{\xb}{\mathbf{x}}
\newcommand{\yb}{\mathbf{y}}
\newcommand{\zb}{\mathbf{z}}

\newcommand{\eb}{\mathbf{e}}

\newcommand{\Hb}{\mathbf{H}}

\newcommand{\Ab}{\mathbf{A}}
\newcommand{\Bb}{\mathbf{B}}
\newcommand{\Cb}{\mathbf{C}}
\newcommand{\Lb}{\mathbf{L}}

\newcommand{\Wb}{\mathbf{W}}

\newcommand{\Ub}{\mathbf{U}}
\newcommand{\Vb}{\mathbf{V}}
\newcommand{\Pb}{\mathbf{P}}
\newcommand{\Qb}{\mathbf{Q}}
\newcommand{\Ib}{\mathbf{I}}
\newcommand{\Mb}{\mathbf{M}}

\newtheorem{remark}{Remark}

\newtheorem{lemma}{Lemma}
\newtheorem{proposition}{Proposition}

\newtheorem{assumption}{Assumption}

\newtheorem{step}{Step}

\allowdisplaybreaks

\begin{document}

\begin{frontmatter}

\title{Distributed Least Squares Solver for Network Linear Equations}
\vspace{-0.4cm}%

\thanks[CorrespondingAuthor]{Corresponding author. Tel.: +86-551-63601508.}

\vspace{-0.4cm}
\author[NEU]{Tao Yang}\ead{yangtao@mail.neu.edu.cn},
\author[ARL]{Jemin George}\ead{jemin.george.civ@mail.mil},
\author[USTC]{Jiahu Qin\thanksref{CorrespondingAuthor}}\ead{jhqin@ustc.edu.cn},
\author[KTH RIT]{Xinlei Yi}\ead{xinleiy@kth.se},
\author[ZJU]{Junfeng Wu}\ead{jfwu@zju.edu.cn}

\vspace{-0.2cm}
\address[NEU]{State Key Laboratory of Synthetical Automation for Process Industries, Northeastern University, Shenyang 110819, China}
\vspace{-0.2cm}
\address[ARL]{U.S. Army Research Laboratory, Adelphi, MD 20783, USA}
\vspace{-0.2cm}
\address[USTC]{Department of Automation, University of Science and Technology of China, Hefei 230027, China}
\vspace{-0.2cm}
\address[KTH RIT]{Division of Decision and Control Systems, School of Electrical Engineering and Computer Science, KTH
	Royal Institute of Technology, 100 44, Stockholm, Sweden}
\vspace{-0.2cm}
\address[ZJU]{College of Control Science and Engineering, Zhejiang University, Hangzhou 310027, China}

\vspace{-0.8cm}

\begin{abstract}
In this paper, we study the problem of finding the least square solutions of over-determined linear algebraic equations over networks in a distributed manner.
Each node has access to one of the linear equations and holds a dynamic state.
We first propose a distributed least square solver over connected undirected interaction graphs and establish a necessary and sufficient on the step-size under which the algorithm exponentially converges to the least square solution.
Next, we develop a distributed least square solver over strongly connected directed graphs and show that the proposed algorithm exponentially converges to the least square solution provided the step-size is sufficiently small. 
Moreover, we develop a finite-time least square solver by equipping the proposed algorithms with a finite-time decentralized computation mechanism.
The theoretical findings are validated and illustrated by numerical simulation examples.
\vspace{-0.5cm}
\end{abstract}
\end{frontmatter}

\section{Introduction}\label{sec-intro}
\vspace{-3mm}
\setlength{\parindent}{3ex}
In recent years, the development of distributed algorithms to solve linear algebraic equations over multi-agent networks has attracted much attention due to the fact that linear algebraic equations are fundamental for various practical engineering applications \citep{Mou-TAC15,Mou_SCL16,anderson2015decentralized,JiLiu-Auto,JiLiu-TAC,Guodong-TAC17,Jie-TCNS2018,PengWang-TCNS,Hong-TAC18}.
In these algorithms, each node has access to one equation and holds a dynamic state, which is an estimate of the solution.
By exchanging their states with neighboring nodes over an underlying interaction graph, all nodes collaboratively solve the linear equations.
Various distributed algorithms based on distributed control and optimization have been developed for solving the linear equations which have exact solutions, among which discrete-time algorithms are given in \citet{Mou-TAC15,JiLiu-TAC,JiLiu-Auto,Jie-TCNS2018,PengWang-TCNS} and continuous-time algorithms are presented in \citet{anderson2015decentralized,Guodong-TAC17}.
However, most of these existing algorithms can only produce least square solutions for over-determined linear equations in the approximate sense \citep{Mou-TAC15} or for limited graph structures \citep{Elia_ACC2012,Guodong-TAC17}.

\vspace{-3mm}
\setlength{\parindent}{3ex}
By reformulating the least square problem as a distributed optimization problem, various distributed optimization algorithms have be proposed.
For example, a continuous-time version of distributed algorithms proposed in \citet{Nedic09,Nedic10} has been applied to solve the exact least square problem in \citet{Guodong-TAC17}.
However, the drawback is the slow convergence rate due to the diminishing step-size.
With a fixed step-size, it can only find approximated least square solutions with a bounded error.
The recent studies focus on developing distributed algorithms with faster convergence rates to find the exact least square solutions, see, e.g., continuous-time algorithms proposed in \citet{Elia_Allerton10,Cortes_TAC_CT,Guodong-arXiv} based on the classical Arrow-Hurwicz-Uzawa flow \citep{Arrow58}, and discrete-time algorithms proposed in \citet{Elia_ACC2012,Guodong-arXiv,Mou-TAC}.

\vspace{-3mm}
\setlength{\parindent}{3ex}
Due to the exponential convergence of these existing algorithms, all nodes need to constantly perform local computation and communicate with their neighboring nodes, which results in a waste of computation and communication resources.
This is not desirable in multi-agent networks since each node is usually equipped with limited communication resources.
Therefore, the fundamental problem is how to find the exact least square solution in a finite number of iterations, and hence terminate further communication and computation to save energy and resources. This motivates our study of this paper.

\vspace{-3mm}
\setlength{\parindent}{3ex}
The contributions of this paper are summarized as follows.
\vspace{-3mm}
\begin{itemize}
\item
First, we develop a distributed algorithm for solving the least square problem over connected undirected graphs.
We explicitly establish a critical value on the step-size, below which the algorithm exponentially converges to the least square solution, and above which the algorithm diverges.
Our proposed algorithm is discrete-time and readily to be implemented,
while the algorithms proposed in \citet{Elia_ACC2012,Guodong-arXiv} are continuous-time
and require the discretization for the implementation.
Compared to existing studies for distributed optimization for strongly convex and smooth local cost functions \citep{Xu_CDC15,Lina-TCNS,Nedic-arXiv1,Dusan-Unified}, which only establish sufficient conditions on the step-size for the exponential convergence, in this paper, we focus on the case where local cost functions are quadratic and only positive semidefinite (convex) but not positive definite (strongly convex).
Moreover, we establish a necessary and sufficient condition on the step-size for the exponential convergence.

\item
Furthermore, we develop a distributed least square solver over directed graphs and show that the proposed algorithm exponentially converges to the least square solution if the step-size is sufficiently small.
Compared with the existing distributed algorithms for computing the exact least square solutions~\citep{Elia_Allerton10,Cortes_TAC_CT,Elia_ACC2012,Mou-TAC,Guodong-arXiv}, which are only applicable to connected undirected graphs or weight-balanced strongly connected digraphs, our proposed algorithm is applicable to strongly connected directed graphs, which are not necessarily weight-balanced.

\item
Last but not least, we develop a finite-time least square solver by equipping the proposed distributed algorithms with a decentralized computation mechanism based on the finite-time consensus technique proposed in~\citet{Sunda07,Ye-AUTO,Themis-TCNS,YangT16,LishaoYao-ICCA2018}.
The proposed mechanism enables an arbitrarily chosen node to compute the exact least square solution within a finite number of iterations, by using its local successive state values obtained from the underlying distributed algorithm.
With the finite-time computation mechanism, nodes can terminate further communication.
This result is among the first distributed algorithms which compute the exact least square solutions in a finite number of iterations.
\end{itemize}

\vspace{-3mm}
\setlength{\parindent}{3ex}
The remainder of the paper is organized as follows:
In Section~\ref{sec-prob}, we formulate the least square problem for linear equations.
In Sections~\ref{sec-convergence}, we present our main results for undirected graphs and directed graphs, respectively.
In Section~\ref{sec-finite-time}, we develop a finite-time least square solver by equipping the proposed algorithms with a decentralized finite-time computation mechanism.
Section \ref{sec-simulation} presents numerical simulation examples.
Finally, concluding remarks are offered in Section~\ref{sec-conclusion}.

\vspace{-3mm}
\section{Problem Formulation}\label{sec-prob}
\vspace{-3mm}
\setlength{\parindent}{3ex}
Consider the following linear algebraic equation with unknown $\yb \in \R^m$:
\begin{equation}\label{linear-eq}
\zb=\Hb \yb,
\end{equation}
where $\zb \in \R^N$ and $\Hb \in \R^{N \times m}$ are known.
It is well known that if ${\bf z} \in \mathrm{span}({\bf H})$, then the linear equation \eqref{linear-eq} always has one or many
exact solutions.
If ${\bf z} \notin \mathrm{span}(\mathbf{H})$, the above equation \eqref{linear-eq} has no exact solution and the least square solution of \eqref{linear-eq} is defined by the solution of the following optimization problem:
\vspace{-2mm}
\begin{equation}\label{LS-OPT}
\min_{{\bf y} \in \R^m}  \frac{1}{2}\| {\bf z}-{\bf H} {\bf y}\|^2.
\end{equation}

\vspace{-2mm}
\begin{assumption}\label{unique-LS-sol-ass}
\vspace{-3mm}
Assume that the matrix $\Hb$ has full column rank, i.e., $\rank(\Hb)=m$.
\end{assumption}

\vspace{-3mm}
\setlength{\parindent}{3ex}
It is well known that under Assumption \ref{unique-LS-sol-ass}, the problem \eqref{LS-OPT} has a unique solution and is given by
\vspace{-2mm}
\begin{equation}\label{LS-sol}
\yb^\ast =({\bf H}^\top {\bf H})^{-1}{\bf H}^\top {\bf z}.
\end{equation}

\vspace{-3mm}
\setlength{\parindent}{3ex}
Denote
\begin{equation}\label{eq-matrixH}
{\bf H}=\begin{bmatrix}
{\bf h}^\top_1 \\
{\bf h}^\top_2 \\
\vdots \\
{\bf h}^\top_N
\end{bmatrix}, \quad {\bf z}=\begin{bmatrix}
z_1 \\
z_2 \\
\vdots \\
z_N
\end{bmatrix},
\end{equation}
where ${\bf h}^\top_i$ is the $i$-th row vector of the matrix ${\bf H}$.
With these notations, we can rewrite the linear equation \eqref{linear-eq} as
\vspace{-2mm}
\[
{\bf h}^\top_i {\bf y}=z_i, \quad i=1,2,\ldots, N.
\]

\vspace{-3mm}
\setlength{\parindent}{3ex}
\vspace{-2mm}
Let $\mathcal{G}=(\mathcal{V},\mathcal{E})$ be an interaction graph with the set of nodes $\mathcal{V}=\{1,2,\dots,N\}$ and the set of edges $\mathcal{E}\subseteq \mathcal{V}\times\mathcal{V}$.
In this paper, we aim to develop a distributed algorithm over the graph $\mathcal{G}$ to compute the least square solution of \eqref{LS-OPT}, where each node $i \in \mathcal{V}$ only has access to the value of ${\hb}_i$ and $z_i$ without knowledge of ${\hb}_j$ or $z_j$ from other nodes.

\vspace{-3mm}
\setlength{\parindent}{3ex}Define ${\bf x}=[{\bf x}_1^\top \ \dots \ {\bf x}_N^\top]^\top \in \R^{Nm}$.
Consider a cost function $F(\cdot): \R^{Nm} \rightarrow \R$:
\vspace{-3mm}
\begin{equation}\label{global-obj}
F({\bf x})=\sum_{i=1}^{N} f_i({\bf x}_i),
\end{equation}
where
\vspace{-1mm}
\begin{equation}\label{local-obj}
f_i({\bf x}_i)= \frac{1}{2}|{\bf h}^\top_i {\bf x}_i-{z}_i|^2.
\end{equation}
The least square problem \eqref{LS-OPT} is equivalent to the following distributed optimization problem
\vspace{-2mm}
\begin{subequations} \label{distributed-opt-view}
\begin{align}
\min_{\xb \in \R^{Nm}} \qquad &  \sum_{i=1}^{N} f_i({\bf x}_i) \\
{\rm s.t.} \qquad  &   \xb_1=\dots=\xb_N.  \label{distributed-opt-view-eq2}
\end{align}
\end{subequations}
Therefore, under Assumption~\ref{unique-LS-sol-ass}, the solution to the above optimization problem is given by $\xb_1=\dots=\xb_N={\yb}^\ast$, where ${\yb}^\ast$ is given by \eqref{LS-sol}.

\vspace{-3mm}
\section{Convergence Results} \label{sec-convergence}
\vspace{-3mm}
\setlength{\parindent}{3ex}
In this section, we solve the least square problem \eqref{LS-OPT} by considering its equivalent problem \eqref{distributed-opt-view} for undirected graphs and directed graphs, respectively.
\vspace{-3mm}
\subsection{Undirected Graphs}\label{sec-undirected}
\vspace{-3mm}
\setlength{\parindent}{3ex}
We first present our proposed algorithm, where each node $i$ maintains two state vectors $\xb_i(t)\in \R^m$ and $\vb_i(t) \in \R^m$, 
which are node $i$'s estimate of the least square solution, and estimate of the average gradient $\frac{1}{N}\sum_{i=1}^{N} \nabla f_i(\xb_i (t))$, respectively. 
At each time step $t$, each node $i$ updates its state vectors as
\vspace{-2mm}
\begin{subequations}\label{distributed-algo}
	\begin{align}
	\xb_i(t+1) =& \sum_{j\in\mathcal{N}_i} W_{ij} \xb_j(t) - \alpha \vb_i(t), \label{distributed-algo-eq1} \\
	\vb_i(t+1) =& \sum_{j\in\mathcal{N}_i} W_{ij} \vb_j(t) \notag \\
& + \nabla f_i(\xb_i (t+1)) -\nabla f_i(\xb_i(t)), \label{PI-algo-eq2}
	\end{align}
\end{subequations}
where the initial condition $\xb_i(0)$ can be chosen arbitrarily and $\vb_i(0)=\nabla f_i(\xb_i(0))$, $\alpha>0$ is the step-size,
$\mathcal{N}_i$ is the set of neighboring nodes of node $i$ including node $i$ itself, i.e.,
$\mathcal{N}_i=\{j| (j,i) \in \mathcal{E} \} \cup \{i\}$,
$W_{ij}$ is the non-negative weight that node $i$ assigns to its neighboring node $j$ in the undirected graph $\mathcal{G}$, and
\vspace{-3mm}
\begin{equation}\label{algo-local-gradient}
\nabla f_i(\xb_i(t))=\hb_i \hb_i^\top \xb_i(t)-z_i \hb_i,
\end{equation}

\vspace{-3mm}
\begin{assumption}\label{undirected-connected-ass}
Assume the interaction graph $\mathcal{G}$ is undirected and connected.
\end{assumption}

\vspace{-3mm}
\begin{assumption}\label{mixing-weight-ass}
The mixing weight matrix $\Wb=[W_{ij}] \in \R^{N \times N}$ satisfies the following properties:
\vspace{-3mm}
\begin{enumerate}[(i)]
\item For any $(i,j) \in \mathcal{E}$, $W_{ij}>0$. Moreover, $W_{ii}>0$ for all $i \in \mathcal{V}$. For other $(i,j)$, $W_{ij}=0$.
\item The matrix $\Wb$ is symmetric and doubly stochastic, that is, $\Wb=\Wb^\top$, $\Wb {\bf 1}_N ={\bf 1}_N$, and ${\bf 1}^\top_N \Wb ={\bf 1}^\top_N$.
\end{enumerate}
\end{assumption}

\vspace{-3mm}
\begin{remark}\label{remark-undirected}
Note that for a connected undirected graph, the matrix $\Wb$ has all its eigenvalues in $(-1,1]$. Moreover, the eigenvalue at $1$ is unique.
Assumption~\ref{mixing-weight-ass} is common in the literature, see, e.g., \citet{Boyd_SCL04,Yin-EXTRA,Lina-TCNS,Nedic-arXiv1,Lihua-Proximal-TAC}.
Several different mixing rules can be used to ensure that all properties of Assumption~\ref{mixing-weight-ass} are satisfied, see \citet{Yin-EXTRA} for details.
For example, one can choose the matrix $\Wb=\mathbf{I}_N-\frac{1}{\tau} \mathbf{L}$, where $\Lb$ is Laplacian matrix associated with the graph $\mathcal{G}$ and $\tau>\frac{1}{2}\lambda_{\max}(\Lb)$,  where $\lambda_{\max}(\mathbf{L})$ is the largest eigenvalue of the Laplacian matrix.
\end{remark}

\vspace{-2mm}
\setlength{\parindent}{3ex}
Let $\xb(t)=[\xb_1^\top(t) \ \dots \ \xb_N^\top(t)]^\top$, $\vb(t)=[\vb_1^\top(t) \ \dots \ \vb_N^\top(t)]^\top$,
$\tilde{\Hb}=\diag(\hb_1\hb_1^\top, \dots, \hb_N\hb_N^\top)$, and
$\zb_H=[z_1\hb_1^\top \ \dots \ z_N \hb_N^\top]^\top$.
Then the algorithm \eqref{distributed-algo} can be rewritten in a compact form: 
\vspace{-2mm}
\begin{subequations}\label{eq:flowModel1}
	\begin{align}
	{\xb}(t+1) &= (\Wb \otimes  \Ib_m) \xb(t) - \alpha \vb(t), \label{eq:flowModel1-eq1} \\
	\vb(t+1) &= (\Wb\otimes \Ib_m) \vb(t)+ \nabla F(\xb(t+1))-\nabla F(\xb(t)),  \label{eq:flowModel1-eq2}
	\end{align}
\end{subequations}
where the initial condition $\vb(0) = \nabla F(\xb(0))$, and
\vspace{-3mm}
\begin{equation}\label{algo-global-grad}
\nabla F(\xb(t))=\tilde{\Hb} \xb(t)-\zb_H.
\end{equation}

\vspace{-5mm}
\begin{remark}\label{remark-undirected-1}
Note that the algorithm \eqref{distributed-algo} or its compact form \eqref{eq:flowModel1} is essentially the same as the algorithms proposed in \citet{Xu_CDC15,Lina-TCNS,Nedic-arXiv1}.
These studies either require all the local cost functions $f_i(\cdot), \, i=1,\ldots, N$ to be strongly convex or at least one local cost function to be strongly convex.
In our case, local cost functions given in \eqref{local-obj} are quadratic, however, they are only convex but not strongly convex since the matrix $\hb_i \hb_i^\top$ is only positive semidefinite but not positive definite.
Therefore, the convergence analysis in these studies cannot be applied. 
\end{remark}

\vspace{-3mm}
\setlength{\parindent}{3ex}
In order to establish the convergence, we rewrite the distributed algorithm \eqref{eq:flowModel1} as the following linear system:
\vspace{-2mm}
\begin{equation}\label{Lina-algo-closed-loop}
\begin{bmatrix}
\xb(t+1) \\
\vb(t+1)
\end{bmatrix}=\Mb \begin{bmatrix}
\xb(t) \\
\vb(t)
\end{bmatrix}, \, \quad \, \vb(0) = \nabla F (\xb(0)),
\end{equation}
where
\vspace{-2mm}
\begin{equation}\label{Lina-closed-matrix}
\Mb =\begin{bmatrix}
\Wb\otimes \Ib_m & -\alpha \Ib_{Nm} \\
-\tilde{\Hb} \big((\Ib_N-\Wb) \otimes \Ib_m\big) & \Wb\otimes \Ib_m-\alpha \tilde{\Hb}
\end{bmatrix}.
\end{equation}

\vspace{-4mm}
\setlength{\parindent}{3ex}
The matrix $\Mb$ has the following property.
\vspace{-3mm}
\begin{lemma}\label{lemma-semisimple-undirected}
Let Assumptions~\ref{unique-LS-sol-ass}--\ref{mixing-weight-ass} hold.
Then regardless of the value of the step-size $\alpha$,  $1$ is an eigenvalue of the matrix $\Mb$ given by~\eqref{Lina-closed-matrix} and its algebraic
multiplicity is equal to the geometric multiplicity, which is $m$. 
\end{lemma}

\vspace{-3mm}
\setlength{\parindent}{3ex}
The proof of Lemma~\ref{lemma-semisimple-undirected} is given in Appendix~\ref{sec-proof-lemma1}.
The convergence of the distributed algorithm~\eqref{eq:flowModel1} (or equivalently the linear system~\eqref{Lina-algo-closed-loop}) depends on the locations of all other non-unity eigenvalues of the matrix $\Mb$, as shown in the following proposition, whose proof is given in Appendix~\ref{sec-proof-undirected-1}.

\vspace{-3mm}
\begin{proposition}\label{thm-undirected-1}
Let Assumptions~\ref{unique-LS-sol-ass}--\ref{mixing-weight-ass} hold.
Then the distributed algorithm \eqref{eq:flowModel1} exponentially converges to the least square solution if and only if the step-size $\alpha$ is selected such that all other non-unity eigenvalues of the matrix $\Mb$ given by~\eqref{Lina-closed-matrix}, except the $m$ semisimple eigenvalues at $1$,
are strictly within the unit circle.
\end{proposition}

\vspace{-3mm}
\setlength{\parindent}{3ex}
The necessary and sufficient condition given in Proposition~\ref{thm-undirected-1} is implicit.
The following theorem whose proof is given in Appendix~\ref{sec-proof-undirected},
establishes an explicit condition on the step-size.

\vspace{-3mm}
\begin{thm}\label{thm-undirected}
Let Assumptions~\ref{unique-LS-sol-ass}--\ref{mixing-weight-ass} hold.
Then the distributed algorithm \eqref{eq:flowModel1} exponentially converges to the least square solution if and only if $\alpha < \bar{\alpha}$, where
\vspace{-2mm}
\begin{equation}\label{critical-value-undirected}
\bar{\alpha}=\frac{1}{2 \lambda_{\max}\Big(\big((\Ib_N+\Wb)^{-2} \otimes \Ib_m \big)\tilde{\Hb} \Big)}.
\end{equation}
\end{thm}

\vspace{-3mm}
\begin{remark}\label{critical-value-remark}
Theorem~\ref{thm-undirected} explicitly characterizes the critical value on the step-size, below which the algorithm exponentially converges to the least square solution, and above which the algorithm diverges.
The explicit critical value depends on the mixing weight matrices $\Wb$ and the matrix $\Hb$.
Note that the existing studies  \citep{Xu_CDC15,Lina-TCNS,Nedic-arXiv1} only established conservative sufficient conditions on the step-size for the exponential convergence.
Theorem~\ref{thm-undirected} provides a necessary and sufficient condition on the step-size for quadratic cost functions.
\end{remark}

\blue{
\begin{remark}
Also note that similar to the existing studies \citep{Xu_CDC15,Lina-TCNS,Nedic-arXiv1}, the upper bound $\bar{\alpha}$ of the step-size cannot be computed exactly in a distributed way. 
Nevertheless, for the case that the mixing weight $\Wb=\Ib_N-\frac{1}{\tau} \Lb$, where $\tau=\max_{i \in \mathcal{V}} \{d_i\}+1$ and $d_i$ is the weighted degree of node $i$, we can estimate a lower bound of the critical value, which in turn provides a sufficient condition for the proposed algorithm ~\eqref{eq:flowModel1} to exponentially converge to the least square solution as follows. 
Note that from \eqref{critical-value-undirected}, we have $\bar{\alpha} \geq \frac{\lambda^2_{\min}(\Ib_N+\Wb)}{2\lambda_{\max}(\tilde{\Hb})}$. 
Next, since $\lambda_{\max}(\Lb) \leq 2 \max_{i \in \mathcal{V}} \{d_i\}$ (see, eq. (12) of \citet{Olfati07}), we have $\Ib_N+\Wb=2 \Ib_N-\frac{1}{\tau} \Lb \geq  \frac{2}{\max_{i \in \mathcal{V}} \{d_i\} +1} \Ib_N$. 
Therefore, $\lambda^2_{\min}(\Ib_N+\Wb)\geq  \frac{4}{(\max_{i \in \mathcal{V}} \{d_i\} +1)^2}$. 
Hence, we find the lower bound of the critical value $\bar{\alpha} \geq \frac{2}{(\max_{i \in \mathcal{V}} \{d_i\} +1)^2 \lambda_{\max}(\tilde{\Hb})}$. 
Note that both $\max_{i \in \mathcal{V}} \{d_i\}$ and $\lambda_{\max}(\tilde{\Hb})$ can be obtained in a distributed manner and in finite-time by using the max-consensus algorithm proposed in \citet{Minimum-consensus}. 
Therefore, we find a more conservative upper bound for the step-size $\alpha$ to ensure the exponential convergence of the proposed algorithm in a distributed manner, that is,
\vspace{-3mm}
\begin{equation}\label{conservative-upper-bound}
\alpha< \frac{2}{(\max_{i \in \mathcal{V}} \{d_i\} +1)^2 \lambda_{\max}(\tilde{\Hb})}.
\end{equation}
\end{remark}
}

\vspace{-2mm}
\subsection{Directed Graphs}\label{sec-directed}
\vspace{-3mm}
\setlength{\parindent}{3ex} In this section, we extend the algorithm \eqref{eq:flowModel1} to handle directed graphs.
Our proposed algorithm is based on recently developed distributed optimization algorithms for directed graphs~\citep{Yang-PESGM2018,Usman-CSL2018,pu2018push}.
Rather than using the doubly stochastic matrix $\Wb$ as in \eqref{eq:flowModel1}, the proposed algorithm
uses a row stochastic matrix for the mixing of estimates of the least square solution in the update \eqref{eq:flowModel1-eq1}, and employs a column stochastic matrix for tracking the average gradient in the update \eqref{eq:flowModel1-eq2}.
More specifically, at time step $t$, each node $i$ performs the following updates:
\vspace{-2mm}
\begin{subequations}\label{algo-push-pull}
\begin{align}
\xb_i(t+1)&= \sum_{j \in \mathcal{N}^{\text{in}}_i}p_{ij} \xb_j(t)-\alpha \vb_i(t), \label{algo-push-pull-1}\\
\vb_i(t+1)&=  \sum_{j \in \mathcal{N}^{\text{in}}_i} q_{ij} \vb_j(t)+ \nabla f_i(\xb_i(t+1)) - \nabla f_i(\xb_i(t)),  \label{algo-push-pull-2}
\end{align}
\end{subequations}
where the initial condition $\xb_i(0)$ can be chosen arbitrarily and $\vb_i(0)=\nabla f_i(\xb_i(0))$, and $\mathcal{N}^{\text{in}}_i=\{ j \in \mathcal{V} \; | \; (j,i) \in \mathcal{E}\} \cup \{i\}$ is the in-neighbor set of node $i$.

\vspace{-3mm}
\begin{assumption}\label{directed-SC-ass}
Assume the interaction graph $\mathcal{G}$ is directed and strongly connected.
\end{assumption}

\vspace{-3mm}
\begin{assumption}\label{mixing-weight-PQ-ass}
The mixing weight matrices $\Pb=[p_{ij}] \in \R^{N \times N}$ and $\Qb=[q_{ij}] \in \R^{N \times N}$ satisfy the following properties:
\vspace{-3mm}
\begin{enumerate}[(i)]
\item $\Pb$ is row stochastic and $\Qb$ is column stochastic.
\item $p_{ij}>0$ if $j \in \mathcal{N}^{\text{in}}_i$, and $p_{ij}=0$ otherwise.
\item $q_{ij}>0$ if $i \in \mathcal{N}^{\text{out}}_j$, where $\mathcal{N}^{\text{out}}_j=\{ i \in \mathcal{V} \; | \; (j,i) \in \mathcal{E}\} \cup \{j\}$ is the out-neighbor set of node $j$, and $q_{ij}=0$ otherwise.
\end{enumerate}
\end{assumption}

\vspace{-3mm}
\setlength{\parindent}{3ex}
Several choices of the weight matrices $\Pb$ and $\Qb$ which satisfy Assumption~\ref{mixing-weight-PQ-ass} are discussed in \citet{Yang-PESGM2018,Usman-CSL2018,pu2018push}.
One particular choice is
\vspace{-2.5mm}
\begin{equation}\label{edge-weight-PQ}
p_{ij}=\begin{cases}
\frac{1}{d^{\text{in}}_i} & \mbox{if } j \in \mathcal{N}^{\text{in}}_i \\
0 & \mbox{otherwise }
\end{cases}, \, \,
q_{ij}=\begin{cases}
\frac{1}{d^{\text{out}}_j} & \mbox{if } i \in \mathcal{N}^{\text{out}}_j  \\
0 & \mbox{otherwise }
\end{cases},
\end{equation}
where $d^{\text{in}}_i$ and $d^{\text{out}}_j$ are in-degree and out-degree of node $i$ and node $j$, respectively.

\vspace{-2mm}
\setlength{\parindent}{3ex}
Note that algorithm \eqref{algo-push-pull} can be written in a compact form as
\vspace{-3mm}
\begin{subequations}\label{eq:flowModel4}
\begin{align}
\xb(t+1) &= (\Pb \otimes  \Ib_m) \xb(t) - \alpha \vb(t), \label{eq:flowModel4-eq1} \\
\vb(t+1) &= (\Qb \otimes \Ib_m) \vb(t)+ \nabla F(\xb(t+1))-\nabla F(\xb(t)), \label{eq:flowModel4-eq2}
\end{align}
\end{subequations}
where $\vb(0) = \nabla F(\xb(0))$. 
Also note that the distributed algorithm \eqref{eq:flowModel4} can be written as the linear system of the form ~\eqref{Lina-algo-closed-loop}, however with a different system matrix:
\vspace{-2mm}
\begin{equation}\label{Lina-closed-matrix-directed}
\Mb =\begin{bmatrix}
\Pb\otimes \Ib_m & -\alpha \Ib_{Nm} \\
-\tilde{\Hb} \big((\Ib_N-\Pb) \otimes \Ib_m\big) & \Qb\otimes \Ib_m-\alpha \tilde{\Hb}
\end{bmatrix}.
\end{equation}

\vspace{-3mm}
\setlength{\parindent}{3ex}
The following lemma whose proof is given in Appendix~\ref{sec-proof-lemma2}, shows that regardless of the value of the step-size $\alpha>0$, the matrix $\Mb$ given in \eqref{Lina-closed-matrix-directed} only has $m$ semisimple eigenvalues at $1$.

\vspace{-3mm}
\begin{lemma}\label{lemma-semisimple-directed}
Let Assumptions~\ref{unique-LS-sol-ass},~\ref{directed-SC-ass} and~\ref{mixing-weight-PQ-ass} hold.
Then regardless of the value of the step-size $\alpha$, $1$ is an eigenvalue of the matrix $\Mb$ given by~\eqref{Lina-closed-matrix-directed} and its algebraic
multiplicity is equal to the geometric multiplicity, which is $m$. 
\end{lemma}

\vspace{-3mm}
\setlength{\parindent}{3ex}
The convergence of the distributed algorithm~\eqref{eq:flowModel4} depends on the locations of all other non-unity eigenvalues of the matrix $\Mb$ given by \eqref{Lina-closed-matrix-directed}, as shown in the following proposition, whose proof is similar to the proof of Proposition~\ref{thm-undirected-1}, and thus omitted.

\vspace{-2mm}
\begin{proposition}\label{directed-thm1}
Let Assumptions~\ref{unique-LS-sol-ass},~\ref{directed-SC-ass} and~\ref{mixing-weight-PQ-ass} hold.
Then the distributed algorithm \eqref{eq:flowModel4} exponentially converges to the least square solution if and only if the step-size $\alpha$ is selected such that all other non-unity eigenvalues of the matrix $\Mb$ given by~\eqref{Lina-closed-matrix-directed}, except the $m$ semisimple eigenvalues at $1$, are strictly within the unit circle.
\end{proposition}

\vspace{-3mm}
\begin{remark}
Note that the algorithm \eqref{eq:flowModel4} is essentially the same as the algorithms proposed in \citet{Yang-PESGM2018,Usman-CSL2018,pu2018push}.
However, the convergence analysis in these studies for general cost functions cannot be applied here due to the same reason as discussed in Remark \ref{remark-undirected-1}.
Compared with these existing studies which established sufficient conditions for the exponential convergence, Proposition~\ref{directed-thm1} provides  a necessary and sufficient conditions for quadratic cost functions.
\end{remark}

\vspace{-3mm}
\setlength{\parindent}{3ex}
The necessary and sufficient condition given in Proposition~\ref{directed-thm1} is implicit.
Unlike the case of undirected graphs, the explicit characterization of a necessary and sufficient condition on the step-size for directed graphs is rather challenging and will be investigated in the future.
Nevertheless, the following theorem whose proof is given in Appendix~\ref{sec-proof-directed}, shows that the algorithm exponentially converges to the least square solution if the step-size is sufficiently small.

\vspace{-3mm}
\begin{thm}\label{directed-thm2}
Let Assumptions~\ref{unique-LS-sol-ass},~\ref{directed-SC-ass} and~\ref{mixing-weight-PQ-ass} hold.
Then the distributed algorithm \eqref{eq:flowModel4} exponentially converges to the least square solution if $\alpha$ is sufficiently small.
\end{thm}

\vspace{-3mm}
\section{Decentralized Finite-time Computation} \label{sec-finite-time}
\vspace{-3mm}
\setlength{\parindent}{3ex}
In this section, we develop a finite-time least square solver by equipping the algorithm~\eqref{eq:flowModel1} for undirected graphs (or the algorithm~\eqref{eq:flowModel4} for directed graphs)  with a decentralized computation mechanism,  which enables an arbitrarily chosen node to compute the exact least square solution in a finite number of time steps, by using the successive values of its local states.

\vspace{-3mm}
\setlength{\parindent}{3ex}
Consider the distributed algorithm~\eqref{Lina-algo-closed-loop} with the matrix $\Mb$ given by \eqref{Lina-closed-matrix} for undirected graphs or \eqref{Lina-closed-matrix-directed} for directed graphs.
Assume that at time step $t$, an arbitrarily chosen node $r \in \mathcal{V}$ has observations about its state $\xb_r(t) \in \R^m$. That is,
\vspace{-2mm}
\begin{equation}\label{pair11}
\yb_r(t)=\begin{bmatrix}
y_{r,1}(t) \\
\vdots \\
y_{r,m}(t)
\end{bmatrix}=\begin{bmatrix}
x_{r,1}(t) \\
\vdots \\
x_{r,m}(t)
\end{bmatrix}=\Cb_{r} \begin{bmatrix}
\xb(t) \\
\vb(t)
\end{bmatrix},
\end{equation}
with
\begin{equation}\label{pair1}
\Cb_{r}=\begin{bmatrix}
\Cb_{r,1} \\
\vdots \\
\Cb_{r,m}
\end{bmatrix}
=\begin{bmatrix}
{\eb}^\top_r \otimes \Ib_m & {\bf 0}_{m \times Nm}
\end{bmatrix},
\end{equation}
where ${\eb}^\top_r \in {\R}^{1 \times N}$ is the row vector whose $r$-th entry is $1$, and the remaining entries are all zeros.

\vspace{-3mm}
\setlength{\parindent}{3ex}Based on these local successive observations, we will propose a decentralized computation algorithm which enables an arbitrarily chosen node $r$ to compute the exact least square solution $\yb^\ast=[y^\ast_1 \ldots y^\ast_m]^\top \in \R^m$ in a finite number of iterations.
The algorithm is motivated by the finite-time technique originally proposed in \citet{Sunda07,Ye-AUTO} for distributed consensus. 
Here, we extend it for computing the exact least square solution. 
To our best knowledge, the finite-time least square solvers are not available in the literature.

\vspace{-3mm}
\setlength{\parindent}{3ex}In order for an arbitrarily chosen node $r \in \mathcal{V}$ to compute the $j$-th component of the least square solution $y^\ast_j$, node $r$ needs to store successive observations for a few time steps.
Consider the $2k+2$ successive observations $y_{r,j}(t)=x_{r,j}(t)$ at node $r$,
that is, $y_{r,j}(0), \, y_{r,j}(1), \, \ldots, y_{r,j}(2k), \, y_{r,j}(2k+1)$.
With these observations, node $r$ then calculates the difference between successive values of $y_{r,j}(t)=x_{r,j}(t)$ according to
\vspace{-2mm}
\begin{equation}\label{Xbar}
\overline{y}_{r,j}(t) \triangleq y_{r,j}(t)-y_{r,j}(t-1), \, \quad t=1,\ldots, 2k+1,
\end{equation}
and constructs a $(k+1) \times (k+1)$ square Hankel matrix as
\vspace{-2mm}
\begin{equation}\label{Hankel-Xbar}
{\Hb}^{(r)}_{j,2k} \triangleq
\begin{bmatrix}
\overline{y}_{r,j}(1) & \overline{y}_{r,j}(2) & \ldots & \overline{y}_{r,j}(k+1) \\
\overline{y}_{r,j}(2) & \overline{y}_{r,j}(3) & \ldots & \overline{y}_{r,j}(k+2) \\
\vdots & \vdots & \ddots & \vdots \\
\overline{y}_{r,j}(k+1) & \overline{y}_{r,j}(k+2) & \ldots & \overline{y}_{r,j}(2k+1)
\end{bmatrix}.
\end{equation}

\vspace{-3mm}
\setlength{\parindent}{3ex}
\blue{A decentralized finite-time computation mechanism is summarized in Algorithm~\ref{algo-finite-time-optimality}. 
The next theorem shows that this algorithm enables any node to compute the exact least square solution by using its own local successive observations in a finite number of iterations. The proof readily follows from \citet{Ye-AUTO}, and is thus omitted.}

\begin{algorithm}
	\caption{Decentralized Finite-time Computation for the Least Square Solution}
	\label{algo-finite-time-optimality}
	\textbf{Data:} Successive observations of $\yb_r(t)=\xb_r(t),~t \in \Z_+$.\\
	\textbf{Result:} The least square solution $\yb^*$.
	
	For each $j=1,\ldots,m$, do the following steps:
	\begin{step}
		Compute the vector of differences $\overline{y}_{r,j}(t)$ by \eqref{Xbar}.
	\end{step}
	\begin{step}
		Construct the square Hankel matrix ${\Hb}^{(r)}_{j,2k}$ by \eqref{Hankel-Xbar}.
	\end{step}
	\begin{step}
		Increase the dimension $k$ of the square Hankel matrix ${\Hb}^{(r)}_{j,2k}$.
	\end{step}
	\begin{step}
		When the square Hankel matrix is singular,
		compute the kernel $\bm{\beta}^{(r)}_j=[{\beta}^{(r)}_{j,0} \ \ldots \ \beta^{(r)}_{j,D_{r,j}-1} \ 1]^\top$.
	\end{step}
	\begin{step}
		Compute the $j$-th component of the least square solution $y^\ast_j$ according to \eqref{final-optimal}:
\vspace{-3mm}
\begin{equation}
y^\ast_j=\frac{\yb^\top_{j,D_{r,j}}\bm{\beta}^{(r)}_j}{{\bf 1}^\top \bm{\beta}^{(r)}_j}, \label{final-optimal}
\end{equation}
\vspace{-3mm}
where
\vspace{-2mm}
\begin{equation*}
\yb_{j,D_{r,j}}=\begin{bmatrix}
y_{r,j}(0) \;\;\; y_{r,j}(1) \;\;\; \ldots \;\;\; y_{r,j}(D_{r,j})
\end{bmatrix}^\top.
\end{equation*}			
	\end{step}
\end{algorithm}

\vspace{-2mm}
\begin{thm}\label{thm-finite-time}
Suppose that the algorithm~\eqref{eq:flowModel1} for undirected graphs (or the algorithm~\eqref{eq:flowModel4} for directed graphs) is convergent. 
Then Algorithm~\ref{algo-finite-time-optimality} enables an arbitrarily chosen node $r\in\mathcal{V}$ to compute the exact least square solution from a finite  number of consecutive states $\xb_r(t)$ observed over a range of time-steps.
\end{thm}

\vspace{-2mm}
\begin{remark}
Note that Algorithm~\ref{algo-finite-time-optimality} relies on the analysis of the rank of a square Hankel matrix.
As shown in Theorem~\ref{thm-finite-time}, for node $r$ to compute the $j$-th component of the least square solution $y^\ast_j$, the $(k+1) \times (k+1)$ square Hankel matrix is guaranteed to lose rank when $k=D_{r,j}$.
That is, an arbitrarily chosen node $r$ with successive observations over $2(D_{r,j}+1)$ number of time steps is able to compute the $j$-th component.
Moreover, as shown in \citet{Ye-AUTO}, the number $D_{r,j}+1$ is equal to the rank of the observability matrix associated with the matrix pair $(\Mb,\Cb_{r,j})$.
\end{remark}

%

\vspace{-2mm}
\section{Numerical Examples}\label{sec-simulation}
\vspace{-3mm}
\setlength{\parindent}{3ex}In this section, we provide numerical examples to validate and illustrate our results.

\begin{figure}[!t]
	\centering
	\includegraphics[scale=0.35]{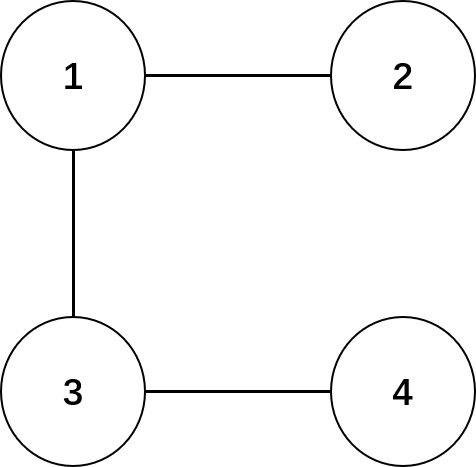}
	\caption{An undirected and connected graph with four nodes.} 
	\label{fig:conv_graph}
\end{figure}

\vspace{-2mm}
\noindent{\bf Example 1.}
In this example, we illustrate the results stated in Theorem~\ref{thm-undirected}.
Consider a linear equation in the form of \eqref{linear-eq} where $\yb\in\R^2$, $\Hb = \bigg[
\begin{smallmatrix}
0\;\; & 1\\
3\;\; & 0\\
2\;\; & 0\\
1\;\; & 0
\end{smallmatrix}\bigg]$ and $\zb = \bigg[
\begin{smallmatrix}
-1\\
0\\
-2\\
2
\end{smallmatrix}\bigg]$.
Since Assumption~\ref{unique-LS-sol-ass} is satisfied, the linear equation has a unique least square solution {$\yb^\ast=[-0.1429 \ -1]^\top$}.
The underlying interaction graph is given in Fig.~\ref{fig:conv_graph}, which is undirected and connected.
Therefore, Assumption~\ref{undirected-connected-ass} is satisfied.
Consider the proposed distributed algorithm \eqref{eq:flowModel1}.
Choose the mixing weight matrix $\Wb= \bigg[
\begin{smallmatrix}
0.7 \;\; & 0.15 \;\; & 0.15 \;\; & 0 \\
0.15\;\;  &  0.85\;\;  & 0\;\;  & 0 \\
0.15\;\;  & 0\;\;  & 0.7\;\;  & 0.15 \\
0\;\;   & 0\;\;  & 0.15\;\;  & 0.85
\end{smallmatrix}\bigg]$.
It is easy to verify that the matrix $\Wb$ satisfies Assumption~\ref{mixing-weight-ass}.
For this case, the critical value given in $\bar{\alpha}$ in \eqref{critical-value-undirected} is $\bar{\alpha}=0.1858$.

\begin{figure}[!t]
	\begin{center}
		\subfigure[$\alpha=0.1857$]
		{
			\includegraphics[scale=0.45]{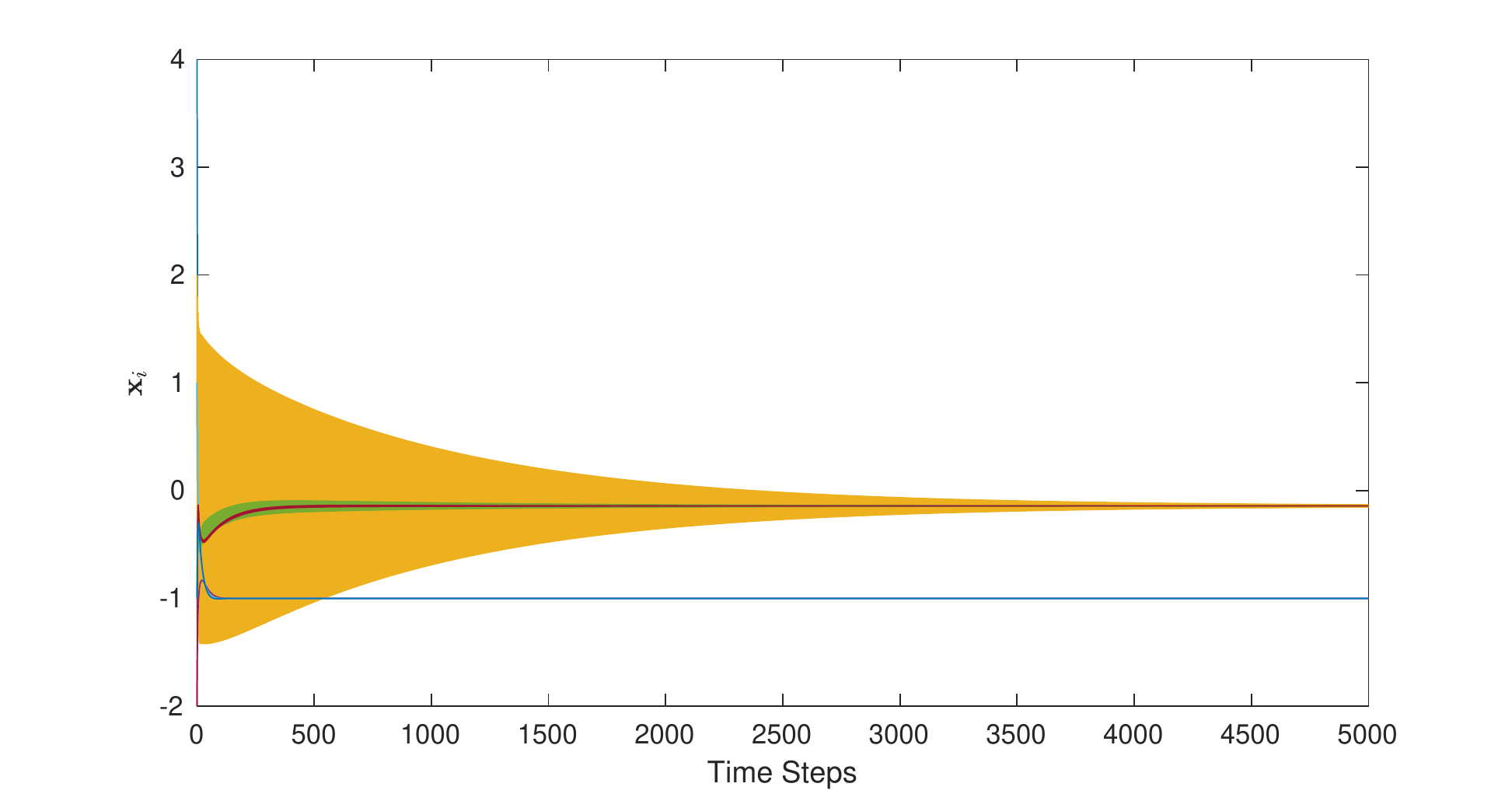}
			\label{fig:our-state-combined-undirected-cov}
		}
		\subfigure[$\alpha=0.1859$]
		{
			\includegraphics[scale=0.45]{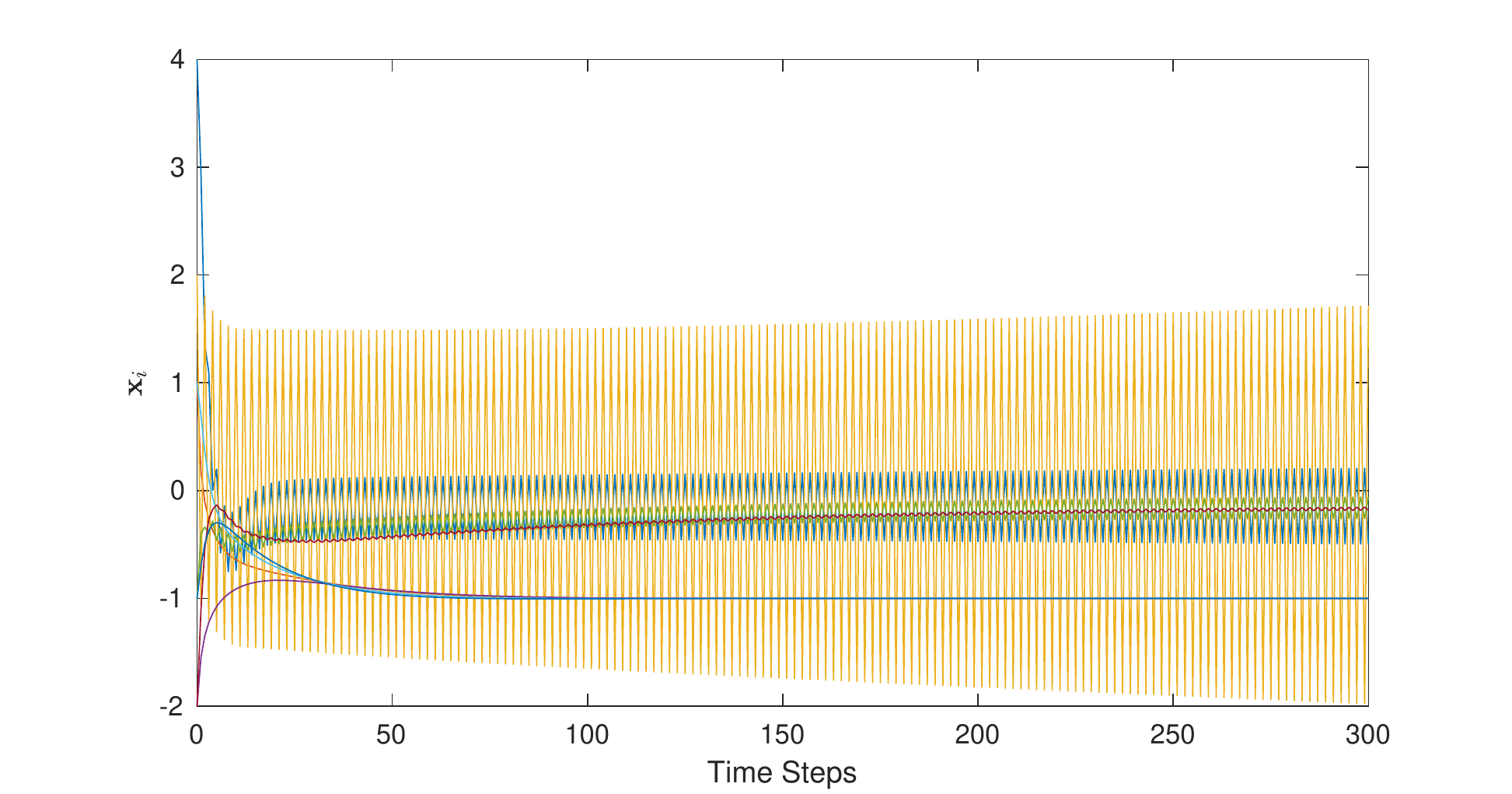}
			\label{fig:our-state-combined-undirected-div}
		}
		\caption{State evolutions $\xb_{i}$ for the algorithm \eqref{eq:flowModel1}.}
	\end{center}
\end{figure}

\vspace{-3mm}
\setlength{\parindent}{3ex}
For the initial condition $\xb(0)=[4 \; 1 \;  2 \;  -2 \;  -1 \;  1 \;  -2 \;  -1]^\top$,
the initial condition $\vb(0)$ is computed as $\vb(0)=\tilde{\Hb}\xb(0)-\zb_H=[0 \;  2 \; 18 \;  0 \;  0 \;  0 \;  -4 \;  0]^\top$.
The state evolutions of $\xb_{i}$ for $i=1,2,3,4$ for the case $\alpha=0.1857$ and $\alpha=0.1859$ are plotted in Fig.~\ref{fig:our-state-combined-undirected-cov} and Fig.~\ref{fig:our-state-combined-undirected-div}, respectively.
It clearly shows that when $\alpha=0.1857<\bar{\alpha}$, all $\xb_i(t)$ converge to the exact least square solution $\yb^\ast$, and when $\alpha=0.1859>\bar{\alpha}$, all $x_{i,1}(t)$ diverge although all $x_{i,2}(t)$ converge to $-1$.
These results are consistent with the results of Theorem~\ref{thm-undirected}.

\begin{figure}[!t]
	\centering
	\subfigure[State evolutions for $\xb_{i}$]
	{
		\begin{picture}(250,120)(0,0)
		\includegraphics[width=0.47\textwidth]{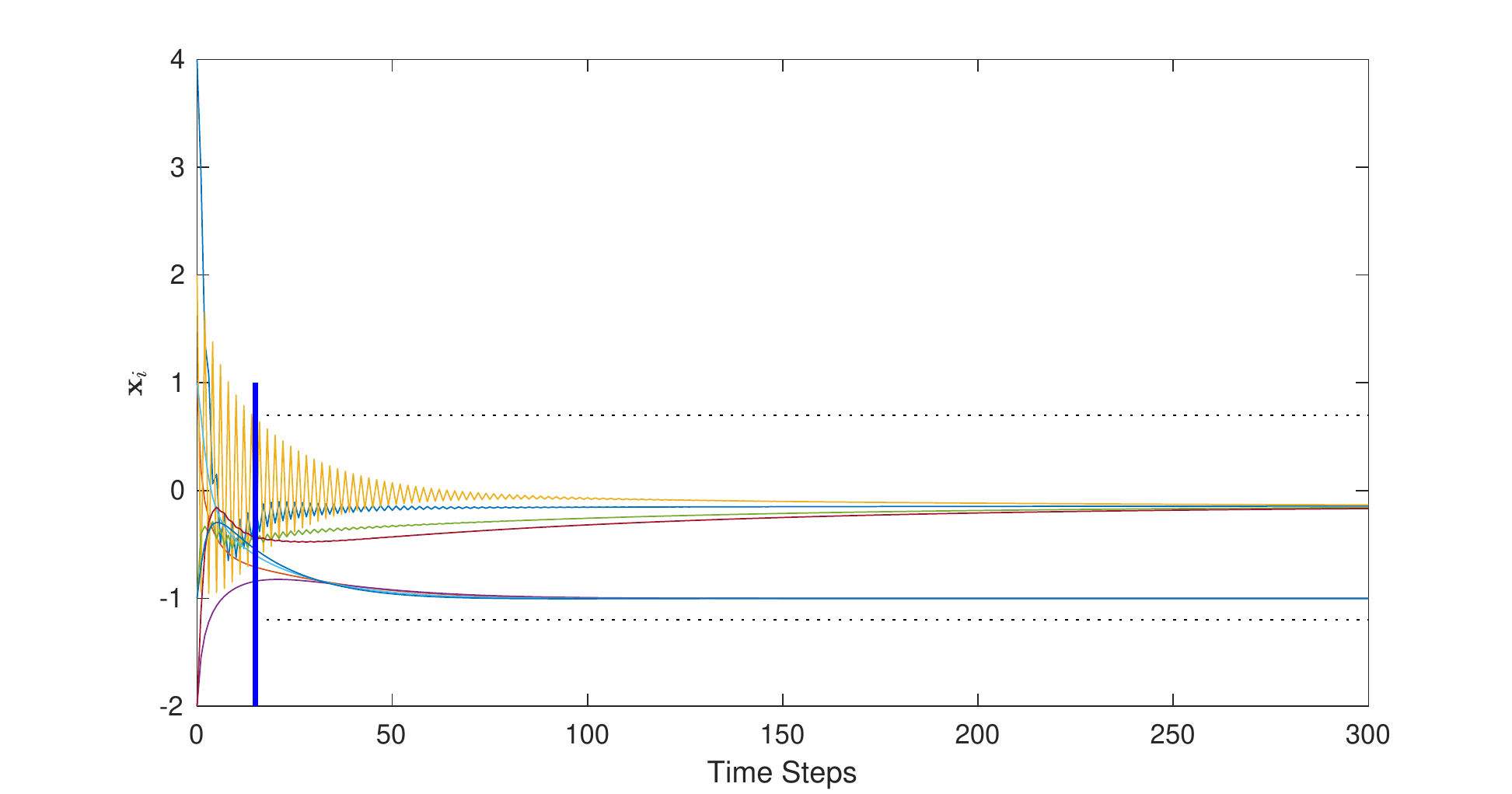}
		\put(-150,60)
		{
			\includegraphics[width=0.20\textwidth]{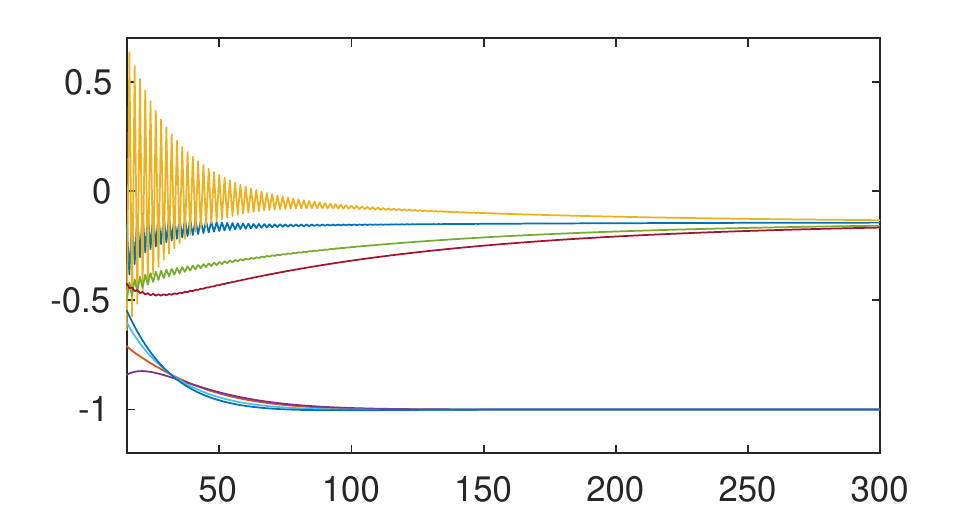}
		}
		\put(-33,59){\vector(-1,1){20}}
		\end{picture}
		\label{fig:FT-combined}
	}
	\subfigure[State evolutions for $x_{i,1}$]
	{
		\begin{picture}(250,120)(0,0)
		\includegraphics[width=0.47\textwidth]{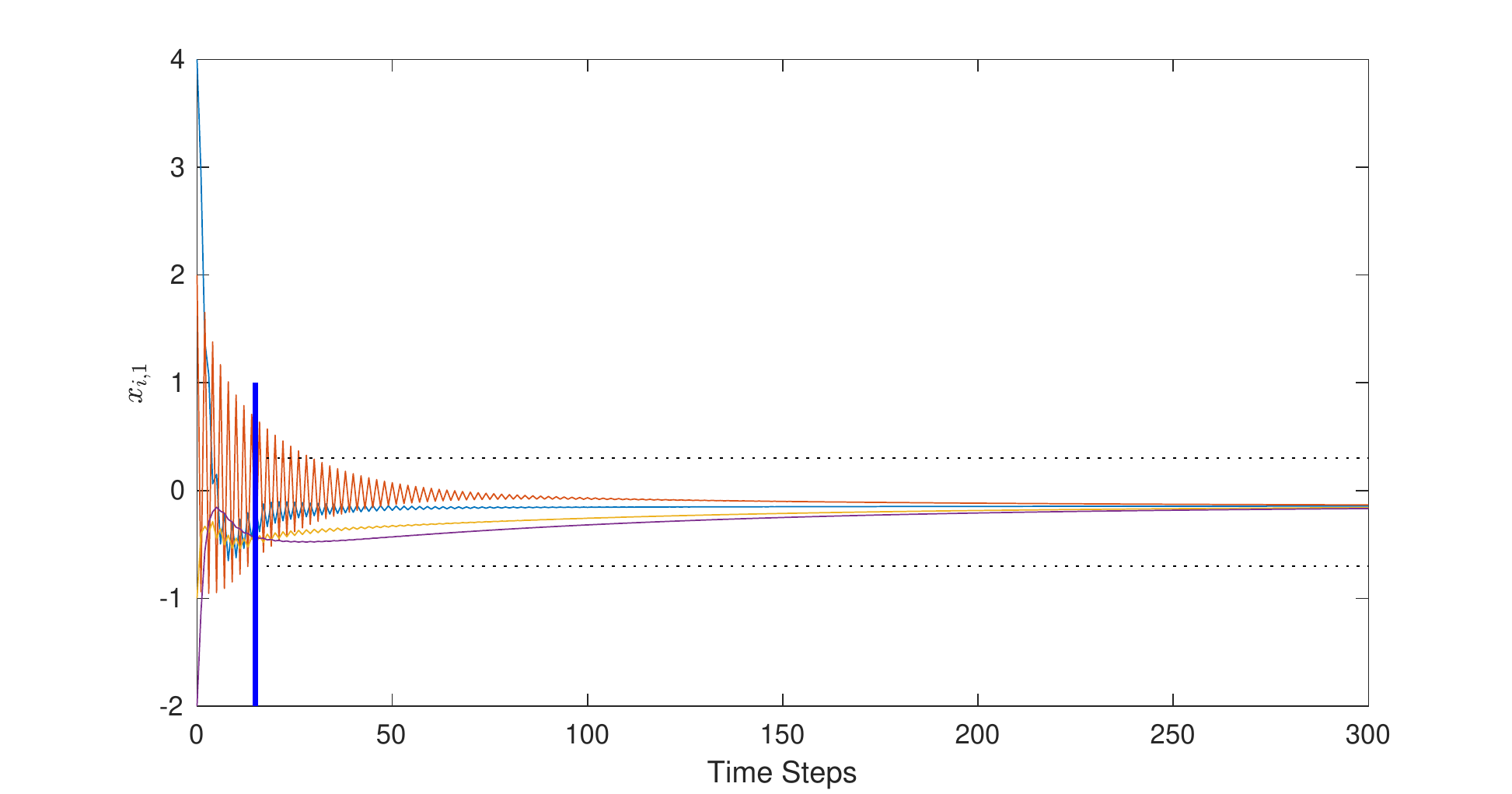}
		\put(-150,60)
		{
			\includegraphics[width=0.20\textwidth]{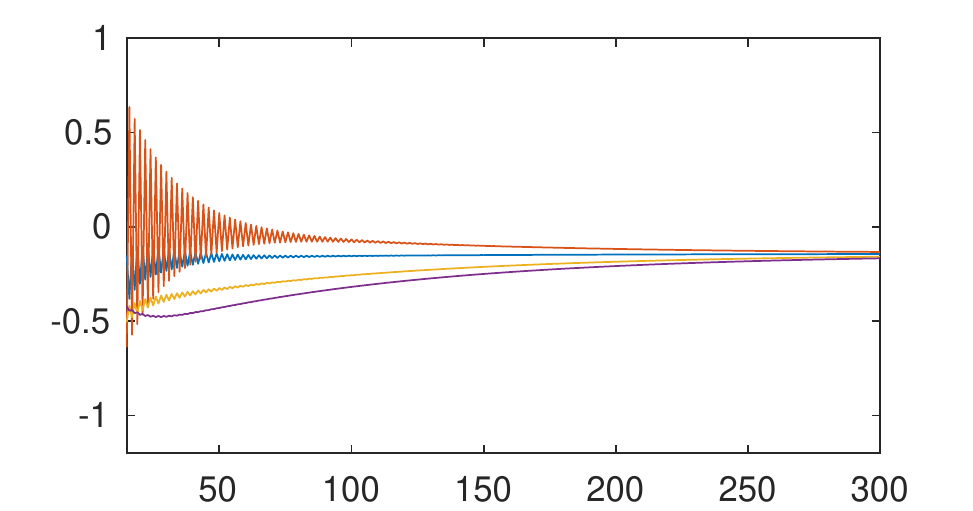}
		}
		\put(-33,59){\vector(-1,1){20}}
		\end{picture}
		\label{fig:FT-x1}
	}
	\caption{Simulation Results. Algorithm \eqref{eq:flowModel1} alone takes roughly $300$ steps to reach the least square solution, while Algorithm \ref{algo-finite-time-optimality} enables nodes to compute the exact least square solution within 16 steps.}
	\label{fig:comparison-FT}
\end{figure}

\vspace{-2mm}
\noindent {\bf Example 2.}
In this example, we illustrate the results stated in Theorem~\ref{thm-finite-time}.
The parameters are chosen the same as those in Example 1.
Now we also equip the algorithm \eqref{eq:flowModel1} with the decentralized finite-time computation mechanism given by Algorithm~\ref{algo-finite-time-optimality}.
The state evolutions of $\xb_i$ for $i=1,2,3,4$ for the case when $\alpha=0.18$ are shown in Fig.~\ref{fig:FT-combined}.
By equipping the algorithm \eqref{eq:flowModel1} with the finite-time mechanism proposed in Algorithm~\ref{algo-finite-time-optimality}, all nodes compute the exact least square solution within $16$ time steps, which is indicated by the vertical blue line in Fig.~\ref{fig:comparison-FT}, and hence terminate further computation and communication.
However, we observe that, at this time step, even approximated least square solution is not achieved by running the algorithm \eqref{eq:flowModel1} alone, and it will take much larger time steps to converge to the least square solution with a reasonable accuracy.


\begin{figure}[!t]
	\centering
	\includegraphics[scale=0.35]{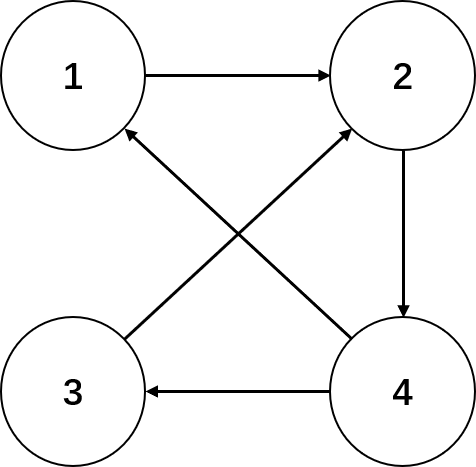}
	\caption{A strongly connected directed graph with four nodes.}
	\label{fig:directed_graph}
\end{figure}

\begin{figure}[!t]
	\centering
	\includegraphics[width=0.47\textwidth]{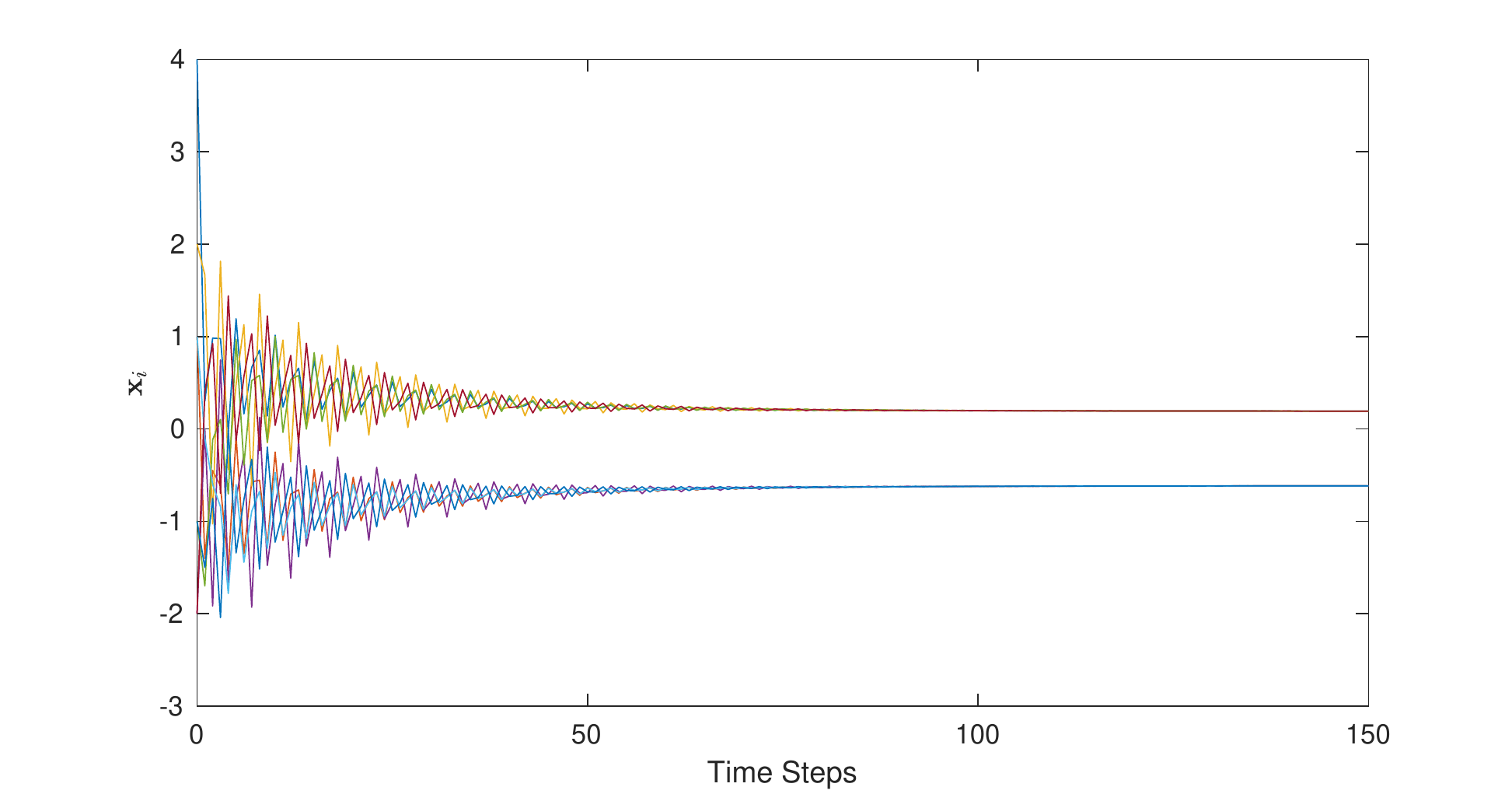}
	\caption{State evolutions of $\xb_{i}$ for algorithm~\eqref{eq:flowModel4} with the step-size $\alpha=0.1$.} 
	\label{fig:our-state-combined-directed}
\end{figure}

\vspace{-3mm}
\noindent{\bf Example 3.}
In this example, we illustrate the results stated in Proposition~\ref{directed-thm1}. 
Consider a linear equation in the form of \eqref{linear-eq} where $\yb\in\R^2$, $\Hb =  \bigg[
\begin{smallmatrix}
1\;\; & 2\\
2\;\; & 2\\
2\;\; & 1\\
1\;\; & 0
\end{smallmatrix}\bigg]$, and $\zb = \bigg[
\begin{smallmatrix}
-1\\
0\\
-2\\
2
\end{smallmatrix}\bigg]$.
Since Assumption~\ref{unique-LS-sol-ass} is satisfied, the linear equation has a unique least square solution $\yb^\ast=[0.1923 \ -0.6514]^\top$.
The directed interaction graph is given in Fig.~\ref{fig:directed_graph}, which is strongly connected. 
Therefore, Assumption~\ref{directed-SC-ass} is satisfied. 
Also note that the digraph is not weight-balanced since $d^{\text{in}}_2=3 \neq d^{\text{out}}_2=2$.
Consider the algorithm \eqref{eq:flowModel4} with the step-size $\alpha=0.1$.
Choose the mixing weight matrices according to \eqref{edge-weight-PQ} as
$\Pb=\Bigg[\begin{smallmatrix}
\frac{1}{2}\;\; & 0\;\; & 0\;\; & \frac{1}{2} \\
\frac{1}{3}\;\; & \frac{1}{3}\;\; & \frac{1}{3}\;\; & 0 \\
0\;\; & 0\;\; & \frac{1}{2}\;\; & \frac{1}{2} \\
0\;\; & \frac{1}{2}\;\; & 0\;\; & \frac{1}{2}
\end{smallmatrix}\Bigg]$ and $\Qb=\Bigg[\begin{smallmatrix}
\frac{1}{2}\;\; & 0\;\;& 0\;\; & \frac{1}{3} \\
\frac{1}{2}\;\; & \frac{1}{2}\;\; & \frac{1}{2} \;\; & 0 \\
0\;\; & 0\;\; & \frac{1}{2}\;\; & \frac{1}{3} \\
0\;\; & \frac{1}{2}\;\; & 0\;\; & \frac{1}{3}
\end{smallmatrix}\Bigg]$,
which satisfy Assumption~\ref{mixing-weight-PQ-ass}.
It is easy to check that all the eigenvalues of the matrix $\Mb$ given by \eqref{Lina-closed-matrix-directed}, except $2$ semisimple eigenvalues at $1$, are strictly within the unit circle. 
Thus according to Proposition~\ref{directed-thm1}, the algorithm \eqref{eq:flowModel4} exponentially converges to the exact least square solution.
The state evolutions of $\xb_{i}$ for $i=1,2,3,4$ are plotted in Fig.~\ref{fig:our-state-combined-directed}, which shows that all $\xb_i(t)$ converge to the exact least square solution $\yb^\ast$.

\begin{figure}[!t]
	\centering
	\subfigure[State evolutions for $x_{i,1}$]
	{
		\begin{picture}(250,120)(0,0)
		\includegraphics[width=0.47\textwidth]{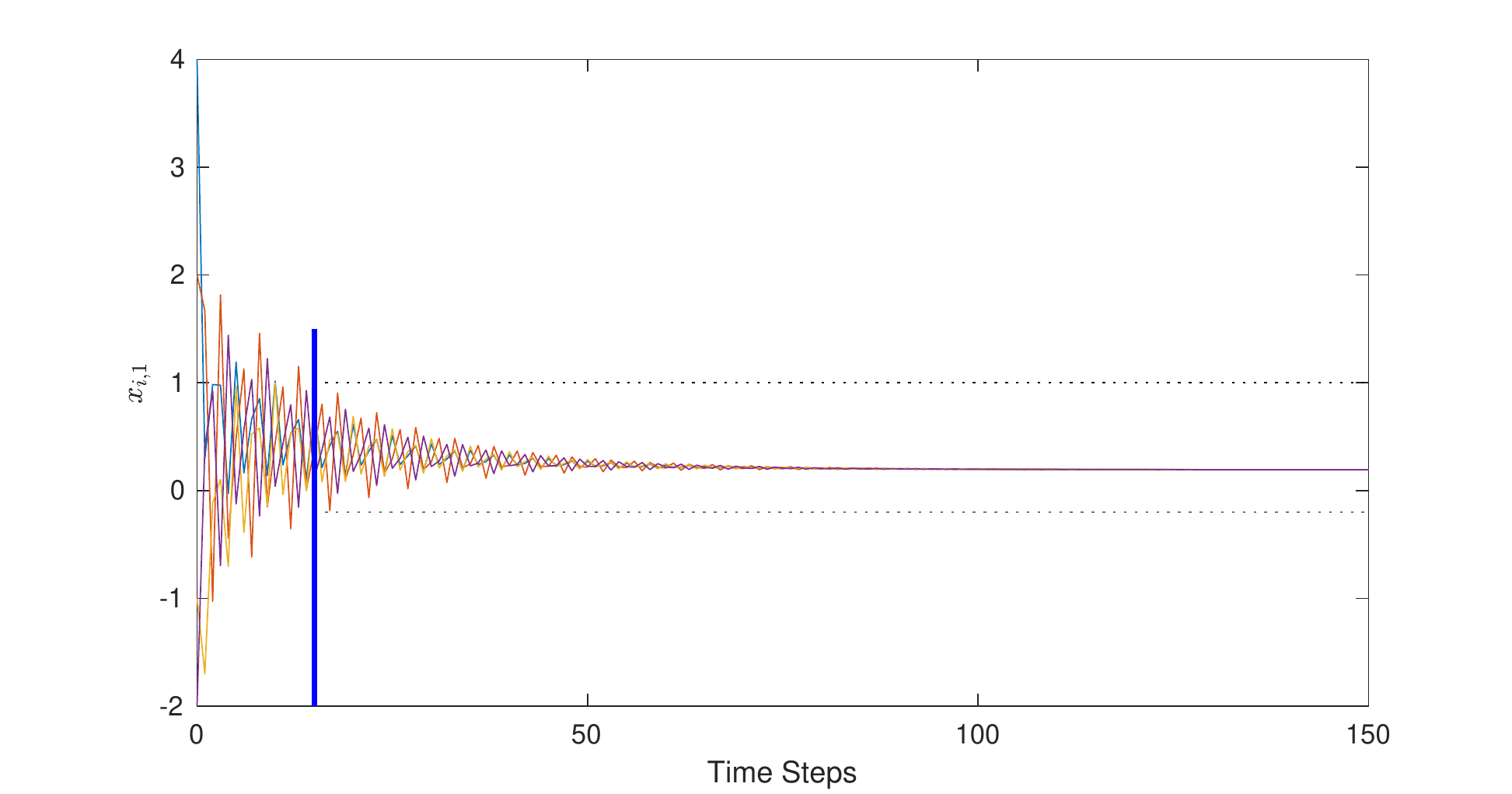}
		\put(-150,60)
		{
			\includegraphics[width=0.20\textwidth]{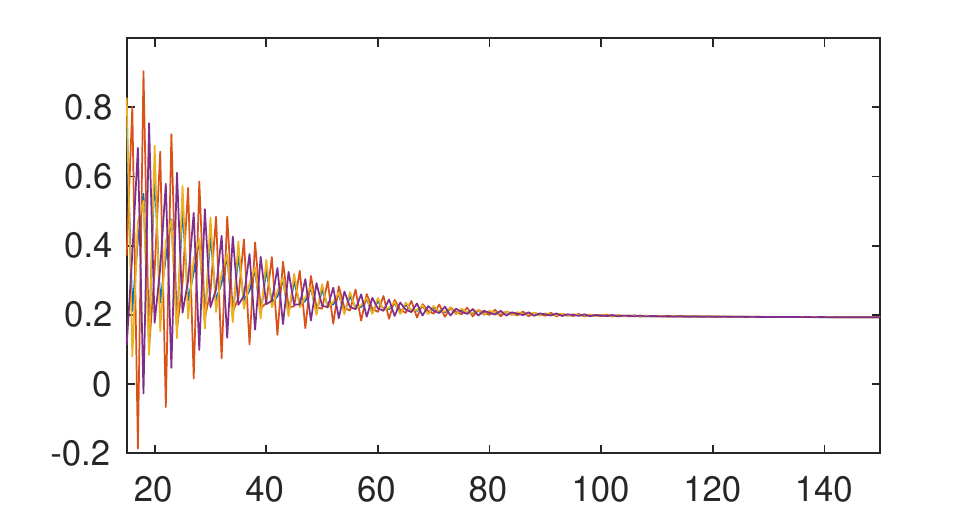}
		}
		\put(-33,59){\vector(-1,1){20}}
		\end{picture}
		\label{fig:example3-x1}
	}
	\subfigure[State evolutions for $x_{i,2}$]
	{
		\begin{picture}(250,120)(0,0)
		\includegraphics[width=0.47\textwidth]{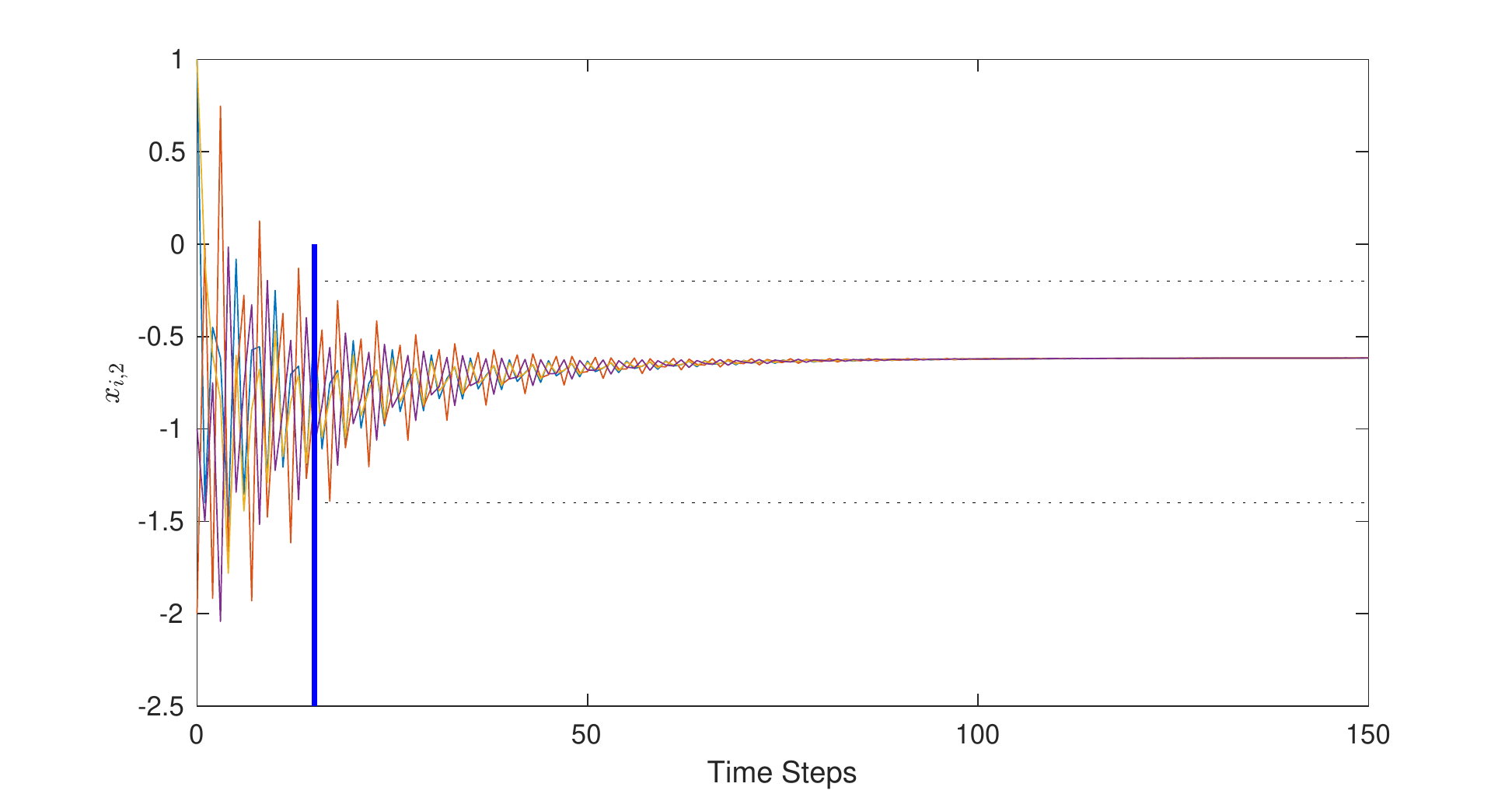}
		\put(-150,60)
		{
			\includegraphics[width=0.20\textwidth]{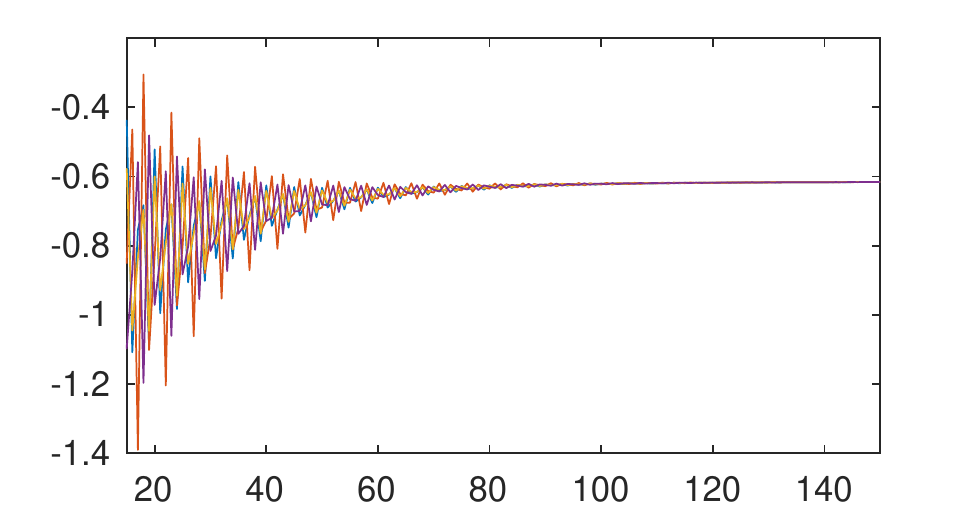}
		}
		\put(-33,59){\vector(-1,1){20}}
		\end{picture}
		\label{fig:example3-x2}
	}
	\caption{Simulation Results. Algorithm \eqref{eq:flowModel4} alone takes roughly $150$ steps to reach the least square solution, while Algorithm \ref{algo-finite-time-optimality} enables nodes to compute the exact least square solution within 16 steps.}
	\label{fig:example-directed-conv}
\end{figure}

\vspace{-3mm}
\setlength{\parindent}{3ex}
Moreover, by equipping the algorithm \eqref{eq:flowModel4} with the finite-time computation mechanism proposed in Algorithm~\ref{algo-finite-time-optimality}, all nodes compute the exact least square solution within $16$ time steps, which are indicated by the vertical blue line in the plots of Fig.~\ref{fig:example-directed-conv}. 
However, we observe that, at this time step, even approximated least square solution is not achieved by running the algorithm \eqref{eq:flowModel4} alone.


\vspace{-3mm}
\section{Conclusions} \label{sec-conclusion}
\vspace{-3mm}
\setlength{\parindent}{3ex}
In this paper, we studied the problem of distributed computing the exact least square solution of over-determined linear algebraic equations over multi-agent networks. 
We first proposed a distributed algorithm as an exact least square solver for undirected connected graphs. 
We established a necessary and sufficient condition on the step-size under which the proposed algorithm exponentially converges to the exact least square solution.
Next, we developed a distributed least square solver for strongly connected directed graphs, which are not necessarily weight-balanced.
We showed that the proposed algorithm exponentially converges to the least square solution if the step-size is sufficiently small. 
Finally, we developed a finite-time exact least square solver for linear equations, by equipping the proposed algorithms with a decentralized computation mechanism.
With the proposed mechanism, an arbitrarily chosen node is able to compute the exact least square solution within a finite number of time steps, by using its local successive observations.
The future direction is to extend the proposed distributed algorithms to networks with time-delays.

\vspace{-2mm}


\appendix
\section{Useful Lemmas}
\vspace{-3mm}
\begin{lemma}\label{lemmaH}
Let Assumption~\ref{unique-LS-sol-ass} holds and $\cb$ be a vector in $\R^m$. If $({\bf 1}_N \otimes \cb)^\top\tilde{\Hb}({\bf 1}_N \otimes \cb)=0$, then $\cb={\bf 0}_{m}$.
\end{lemma}
\vspace{-3mm}
\noindent \textbf{Proof:} Define $\bar{\Hb} \triangleq \diag(\hb^\top_1;\ldots;\hb^\top_N)$.
It is easy to see that $\tilde{\Hb}=\bar{\Hb}^\top \bar{\Hb}$.
Therefore,  $({\bf 1}_N \otimes \cb)^\top\tilde{\Hb}({\bf 1}_N \otimes \cb)=0$ implies $\bar{\Hb} ({\bf 1}_N \otimes \cb)={\bf 0}_{Nm}$, or equivalently $\Hb \cb={\bf 0}_N$.
It then follows from $\rank(\Hb)=m$ that $\cb={\bf 0}_m$.

\vspace{-3mm}
\begin{lemma}\citep[Theorem~7.7.3]{Horn2001}\label{lemma1}
Let $\Ab$ and $\Bb$ be real and symmetric matrices and suppose that $\Ab$ is positive definite.
If $\Bb$ is positive semidefinite, then $\Ab \geq \Bb$ (respectively, $\Ab > \Bb$) if and only if $\lambda_{\max}(\Ab^{-1} \Bb) \leq 1$ (respectively, $\lambda_{\max}(\Ab^{-1} \Bb) < 1$).
\end{lemma}

\vspace{-3mm}
\section{Proof of Lemma~\ref{lemma-semisimple-undirected}} \label{sec-proof-lemma1}
\vspace{-3mm}
\setlength{\parindent}{3ex}
By substituting $\Pb=\Wb$ and $\Qb=\Wb$ into the matrix $\Mb$ given by~\eqref{Lina-closed-matrix-directed}, we get the matrix $\Mb$ given by~\eqref{Lina-closed-matrix}. 
Thus, Lemma~\ref{lemma-semisimple-undirected} is a special case of Lemma~\ref{lemma-semisimple-directed} whose proof is given in Appendix~\ref{sec-proof-lemma2}.

\vspace{-3mm}
\section{Proof of Proposition~\ref{thm-undirected-1}}\label{sec-proof-undirected-1}
\vspace{-3mm}
\noindent \textbf{Sufficiency:}
Suppose that all the eigenvalues of the matrix $\Mb$, except $m$ semisimple eigenvalues at $1$, are strictly within the unit circle.
This together with Lemma~\ref{lemma-semisimple-undirected} implies that the matrix there exists a non-singular matrix $\Vb \in \R^{Nm \times Nm}$ , such that $\Mb = \Vb \big[\begin{smallmatrix}
\Ib_m & \\
& \Mb_s
\end{smallmatrix}\big] \Vb^{-1}$, 
where the eigenvalues of the matrix $\Mb_s$ are all the non-unity eigenvalues of the matrix $\Mb$, which are strictly within the unit circle,
the columns of the matrix $\Vb$ are right eigenvectors and generalized right eigenvectors of the matrix $\Mb$,
the rows of the matrix $\Vb^{-1}$ are left eigenvectors and generalized left eigenvectors of the matrix $\Mb$.
Moreover, the first $m$ columns of the matrix $\Vb$ and the first $m$ rows of the matrix $\Vb^{-1}$ are given by
\vspace{-3mm}
\begin{align}
\tilde{\vb}&=\begin{bmatrix}
{\bf 1}_N \otimes \Ib_m \\
{\bf 0}_N \otimes \Ib_m
\end{bmatrix}, \notag\\
\tilde{\bm{\omega}}^\top&=(\Hb^\top \Hb)^{-1}
\begin{bmatrix}
({\bf 1}^\top_N \otimes \Ib_m) \tilde{\Hb} && -{\bf 1}^\top_N \otimes \Ib_m
\end{bmatrix}.\notag
\end{align}
Therefore, we have
\vspace{-3mm}
\begin{align*}\label{asym-solution}
\lim_{t \rightarrow \infty}
\begin{bmatrix}
\xb(t) \\
\vb(t)
\end{bmatrix}&=\Vb \begin{bmatrix}
\Ib_m & \\
& \Mb^t_s
\end{bmatrix} \Vb^{-1} \\
&= \tilde{\vb} \, \tilde{\bm{\omega}}^\top \begin{bmatrix}
\xb (0) \\
\vb (0)
\end{bmatrix} \nonumber \\
&= \begin{bmatrix}
{\bf 1}_N \otimes \Ib_m \\
{\bf 0}_N \otimes \Ib_m
\end{bmatrix} (\Hb^\top \Hb)^{-1} \Big(({\bf 1}_N^\top \otimes \Ib_m) \tilde{\Hb}\xb(0) \nonumber\\
& \hspace{0.4cm}-({\bf 1}_N^\top \otimes \Ib_m) (\tilde{\Hb}\xb(0)-\zb_H)\Big) \\
&= \begin{bmatrix}
{\bf 1}_N \otimes \yb^\ast \\
{\bf 0}_{Nm}
\end{bmatrix},
\end{align*}
where the second equality follows from the fact that all the eigenvalues of the matrix $\Mb_s$ are strictly within the unit circle, the third equality follows the initialization given by \eqref{algo-global-grad},
and the last equality follows from \eqref{LS-sol}, and the fact that
$({\bf 1}_N^\top \otimes \Ib_m)\zb_H=\Hb^\top \zb$,
due to \eqref{eq-matrixH} and the definitions of $\tilde{\Hb}$ and $\zb_H$.

\vspace{-3mm}
\setlength{\parindent}{3ex}Therefore, for all $i=1,\dots,N$,  $\xb_i(t) \rightarrow \yb^\ast$ and $\vb_i(t) \rightarrow {\bf 0}_m$ asymptotically as $t \rightarrow \infty$.
Note that for linear systems, asymptotic convergence and exponential convergence are equivalent.
Thus, the algorithm exponentially converges to the least square solution.

\vspace{-3mm}
\noindent \textbf{Necessity:}
Note that if the matrix $\Mb$ has an eigenvalue outside the unit circle, then the distributed algorithm~\eqref{Lina-algo-closed-loop} is unstable.

\vspace{-3mm}
\section{Proof of Theorem~\ref{thm-undirected}}\label{sec-proof-undirected}
\vspace{-3mm}
\setlength{\parindent}{3ex}
In view of Proposition~\ref{thm-undirected-1}, in order to prove the theorem, it is equivalent to show that all the
eigenvalues of the matrix $\Mb$, except $m$ semisimple eigenvalues at $1$, are strictly within the unit circle if and only if $\alpha<\bar{\alpha}$.

\vspace{-3mm}
\noindent \textbf{Sufficiency:}
Suppose that the condition $\alpha<\bar{\alpha}$ is satisfied, we need to show that all the eigenvalues of the matrix $\Mb$, except $m$ semisimple eigenvalues at $1$, are strictly within the unit circle. 
For any non-unity eigenvalue $\lambda$ of the matrix $\Mb$, there is a nonzero right eigenvector $[\vb^\top_1 \; \vb^\top_2]^\top$, such that
\vspace{-3mm}
\[
\begin{bmatrix}
\Wb\otimes \Ib_m & -\alpha \Ib_{Nm} \\
-\tilde{\Hb} \big((\Ib_N-\Wb) \otimes \Ib_m\big) & \Wb\otimes \Ib_m-\alpha \tilde{\Hb}
\end{bmatrix} \begin{bmatrix}
\vb_1 \\
\vb_2
\end{bmatrix}=\lambda \begin{bmatrix}
\vb_1 \\
\vb_2
\end{bmatrix}.
\]
We then have
\vspace{-3mm}
\begin{subequations}\label{eig-structure}
\begin{align}
& -\alpha \vb_2 =(\lambda \Ib_{Nm} - \Wb \otimes \Ib_m) \vb_1, \label{eig-structure-eq1} \\
& -\tilde{\Hb}\big((\Ib_N-\Wb) \otimes \Ib_m\big) \vb_1=\big(\lambda \Ib_{Nm}-(\Wb \otimes \Ib_m -\alpha \tilde{\Hb})\big) \vb_2. \label{eig-structure-eq2}
\end{align}
\end{subequations}
Since $\alpha>0$, substituting \eqref{eig-structure-eq1} into \eqref{eig-structure-eq2} to eliminate $\vb_2$ yields:
\begin{equation}\label{QEP-undirected1}
\Qb(\lambda) \vb_1 ={\bf 0}_{Nm},
\end{equation}
where
\vspace{-3mm}
\begin{equation}\label{Qb}
\Qb(\lambda)=\lambda^2 \Ib_{Nm}+\big(\alpha \tilde{\Hb}-2 \Wb \otimes \Ib_m \big)  \lambda+\Wb^2 \otimes \Ib_m -\alpha \tilde{\Hb}.
\end{equation}
Define
\begin{equation}\label{QEP-Qs}
\begin{aligned}
\Qb_2 =& \Ib_{Nm}, \\
\Qb_1 =& \alpha \tilde{\Hb}-2 \Wb \otimes \Ib_m, \\
\Qb_0 =& \Wb^2 \otimes \Ib_m -\alpha \tilde{\Hb}.
\end{aligned}
\end{equation}
Note that the matrices $\Qb_i$ for $i=0,1,2$ are real and symmetric.
Then \eqref{QEP-undirected1} becomes:
\vspace{-3mm}
\begin{equation}\label{QEP-undirected}
(\lambda^2 \Qb_2+ \lambda \Qb_1+\Qb_0) \vb_1 ={\bf 0}_{Nm}.
\end{equation}
Note that \eqref{QEP-undirected} implies that 
\vspace{-3mm}
\begin{equation}\label{scalar-ploy}
p(\lambda) \triangleq  a_2 \lambda^2 + a_1 \lambda+a_0=0,
\end{equation}
where
\begin{equation}\label{jury-coeff}
a_i= \vb^{*}_1 \Qb_i \vb_1, \; \text{ for } i=0,1,2,
\end{equation}
and $\vb^*_1$ is the conjugate transpose of the vector $\vb_1$.
Note that $a_i$ for $i=0,1,2$ are real and $a_2>0$ since $\vb_1$ is nonzero.

\vspace{-3mm}
\setlength{\parindent}{3ex}
Denote $\Lambda$ as the solution set to \eqref{scalar-ploy}. Then $\lambda\in\Lambda$.
According to Jury stability criterion \citep{jury}, $|\tilde{\lambda}|<1,~\forall \tilde{\lambda}\in\Lambda$ if and only if
\vspace{-3mm}
\begin{subequations}\label{jury-cond}
\begin{align}
a_2+a_1+a_0& >  0, \label{jury-cond1} \\
a_2-a_1+a_0 &> 0, \label{jury-cond2} \\
|a_0|& <  a_2. \label{jury-cond3}
\end{align}
\end{subequations}
These equations together with \eqref{QEP-Qs} and~\eqref{jury-coeff} imply that $|\tilde{\lambda}|<1,~\forall \tilde{\lambda}\in\Lambda$  if and only if
\vspace{-3mm}
\begin{subequations}\label{jury-cond-new}
\begin{align}
\vb^{*}_1 \big((\Wb-\Ib_N)^2 \otimes \Ib_m \big) \vb_1 & >  0, \label{jury-cond1-new} \\
\vb^{*}_1 \big((\Wb+\Ib_N)^2 \otimes \Ib_m -2\alpha \tilde{\Hb} \big)\vb_1 &> 0, \label{jury-cond2-new} \\
\vb^{*}_1 \big((\Wb^2-\Ib_N)\otimes \Ib_m -\alpha \tilde{\Hb}\big) \vb_1 &< 0, \label{jury-cond31-new} \\
\vb^{*}_1 \big((\Wb^2+\Ib_N)\otimes \Ib_m -\alpha \tilde{\Hb}\big) \vb_1 & > 0. \label{jury-cond32-new}
\end{align}
\end{subequations}
In the following, with the help of \eqref{jury-cond-new}, we show that $|\lambda|<1$. Thus, the sufficiency is proved.

\vspace{-3mm}
\setlength{\parindent}{3ex}
We first show that the condition \eqref{jury-cond31-new} holds. From  $\Wb^2 \leq \Ib_N$ and $\tilde{\Hb} \geq 0$, we know that $\vb^{*}_1 ((\Wb^2-\Ib_N)\otimes \Ib_m -\alpha \tilde{\Hb}) \vb_1 \le 0$. Suppose that the condition \eqref{jury-cond31-new} does not hold. Then,
$\vb^{*}_1 \big((\Wb^2 -\Ib_N)\otimes \Ib_m \big) \vb_1=0$ and $\vb^{*}_1 \tilde{\Hb}\vb_1=0$.
Since $\Ib_N-\Wb^2$ is positive semidefinite and $\vb^{*}_1 \big((\Wb^2 -\Ib_N)\otimes \Ib_m \big) \vb_1=0$, we know that $\big((\Ib_N-\Wb^2)^{\frac{1}{2}}\otimes \Ib_m \big) \vb_1={\bf 0}_{Nm}$. Thus, $((\Ib_N-\Wb^2)\otimes \Ib_m ) \vb_1={\bf 0}_{Nm}$.
This together with the fact that $\Ib_N+\Wb$ is positive definite and $\Ib_N-\Wb^2=(\Ib_N+\Wb)(\Ib_N-\Wb)$ implies that $\big((\Ib_N-\Wb)\otimes \Ib_m\big) \vb_1={\bf 0}_{Nm}$.
Therefore, if the condition \eqref{jury-cond31-new} does not hold, then $(\Wb \otimes \Ib_m) \vb_1=\vb_1$ and $\vb^{*}_1 \tilde{\Hb}\vb_1=0$.
Note that it follows from $(\Wb \otimes \Ib_m) \vb_1=\vb_1$ that $\vb_1={\bf 1}_N \otimes \cb$ for some $\cb \in \R^m$. It It then follows from Lemma~\ref{lemmaH} that $\cb={\bf 0}_m$.
Thus, $\vb_1={\bf 0}_{Nm}$ which contradicts with the fact that $\vb_1$ is a nonzero eigenvector.
Therefore, the condition \eqref{jury-cond31-new} holds.

\vspace{-3mm}
\setlength{\parindent}{3ex}
Next, we show that the condition \eqref{jury-cond2-new} holds if $\alpha<\bar{\alpha}$, where $\bar{\alpha}$ is given by  \eqref{critical-value-undirected}.
We show this by proving that the matrix $(\Wb+\Ib_N)^2 \otimes \Ib_m -2\alpha \tilde{\Hb}$ is positive definite if $\alpha<\bar{\alpha}$.
Note that the matrix $(\Wb+\Ib_N)^2$ is positive definite due to the fact $\Wb>-\Ib_N$ and the matrix $2\alpha \tilde{\Hb}$ is positive semidefinite.
It then follows from \citet[Theorem~7.7.3]{Horn2001} (recapped in  Lemma~\ref{lemma1}) that $(\Wb+\Ib_N)^2 \otimes \Ib_m -2\alpha \tilde{\Hb}$ is positive definite if and only if $\lambda_{\max}((\Wb+\Ib_N)^{-2} \otimes \Ib_m) 2 \alpha \tilde{\Hb})<1$, which is equivalent to $\alpha<\bar{\alpha}$.

\vspace{-3mm}
\setlength{\parindent}{3ex}
We then note that the condition \eqref{jury-cond32-new} is automatically satisfied since the condition \eqref{jury-cond2-new} holds and $\vb^{*}_1 ((\Wb-\Ib_N)^2 \otimes \Ib_m ) \vb_1\ge  0$.

\vspace{-3mm}
\setlength{\parindent}{3ex}
Finally, if the condition \eqref{jury-cond1-new} holds, then according to Jury stability criterion, we know that $|\tilde{\lambda}|<1,~\forall \tilde{\lambda}\in\Lambda$. Thus, $|\lambda|<1$.
If the condition \eqref{jury-cond1-new} does not hold, i.e. $(\Wb \otimes \Ib_m) \vb_1=  \vb_1$,
then it follows from \eqref{scalar-ploy} and \eqref{jury-coeff} that
\vspace{-3mm}
\begin{equation}\label{tao-new1}
p(\lambda)=\big(\vb^{*}_1 \vb_1(\lambda-1)+\alpha \vb^{*}_1 \tilde{\Hb} \vb_1\big)(\lambda-1)=0.
\end{equation}
Since $\lambda \neq 1$ and $\vb_1$ is a nonzero vector, we have
\vspace{-3mm}
\begin{equation}\label{tao-new2}
\lambda=1-\frac{\alpha \vb^{*}_1 \tilde{\Hb} \vb_1}{\vb^{*}_1 \vb_1}.
\end{equation}
We show that $\lambda \in (-1,1)$.
From $\alpha>0$, $(\Wb \otimes \Ib_m) \vb_1=  \vb_1$ and \eqref{jury-cond31-new}, we know that $\alpha \vb^{*}_1 \tilde{\Hb} \vb_1>0$. Thus, $\lambda<1$.
Next, we show that $\lambda>-1$.
Note that it follows \eqref{tao-new2} that in order to show $\lambda>-1$, it is equivalent to show that
$\vb^{*}_1 (2 \Ib_{Nm}-\alpha \tilde{\Hb}) \vb_1>0$.
This can be shown as follows:
\vspace{-3mm}
\begin{align}
\vb^{*}_1 (2 \Ib_{Nm}-\alpha \tilde{\Hb}) \vb_1 &= \frac{1}{2} \vb^{*}_1 (4 \Ib_{Nm}-2\alpha \tilde{\Hb}) \vb_1 \notag \\
& \hspace{-2cm} \geq  \frac{1}{2}  \vb^{*}_1 \big((\Wb+\Ib_N)^2 \otimes \Ib_m -2\alpha \tilde{\Hb} \big)  \vb_1>0, \notag
\end{align}
where the first inequality follows from the fact that $-\Ib_N<\Wb \leq \Ib_N$ and the second inequality follows from \eqref{jury-cond2-new}.

\vspace{-3mm}
\noindent \textbf{Necessity:}
In order to show the necessity, we need to show that if all the eigenvalues of the matrix $\Mb$, except $m$ semisimple eigenvalues at $1$, are strictly within the unit circle, then $\alpha<\bar{\alpha}$.
We show this by contraposition, that is, if $\alpha \geq\bar{\alpha}$, then apart from the eigenvalues at $1$, the matrix $\Mb$ has at least one eigenvalue either on or outside the unit circle.

\vspace{-3mm}
\setlength{\parindent}{3ex}
To show this, we first note that when $\alpha=\bar{\alpha}$, it follows from  Lemma~\ref{lemma1} that the matrix $(\Wb+\Ib_N)^2 \otimes \Ib_m -2\bar{\alpha} \tilde{\Hb}$ is only positive semidefinite but not positive definite. Thus
$\det\big((\Wb+\Ib_N)^2 \otimes \Ib_m -2\bar{\alpha} \tilde{\Hb}\big)=0$.
Next, we show that when $\alpha=\bar{\alpha}$, the matrix $\Mb$ has an eigenvalue at $-1$.
Since the matrices $\lambda \Ib_{Nm}-\Wb\otimes \Ib_m$ and $\alpha \Ib_{Nm}$ commute, it follows from \citet[Theorem~3]{Silvester-det} that
\vspace{-3mm}
\begin{align}\label{det-M}
& \hspace{0.4cm} \det(\lambda \Ib_{2Nm}-\Mb) \notag \\
& =\det\Big(\big(\lambda \Ib_{Nm}-(\Wb \otimes \Ib_m-\alpha \tilde{\Hb})\big)\big(\lambda \Ib_{Nm}-\Wb \otimes \Ib_m \big) \notag \\
& \hspace{0.8cm}- \alpha \tilde{\Hb} \big((\Ib_N-\Wb) \otimes \Ib_m\big) \Big) \notag \\
&=\det \big(\Qb(\lambda)\big),
\end{align}
where $\Qb(\lambda)$ is given by \eqref{Qb}.
This together with the fact that $\det\big((\Wb+\Ib_N)^2 \otimes \Ib_m -2\bar{\alpha} \tilde{\Hb}\big)=0$ implies that
\vspace{-3mm}
\[
\det(-\Ib_{2Nm}-\Mb) =\det\Big( (\Wb+\Ib_N)^2 \otimes \Ib_m -2\bar{\alpha} \tilde{\Hb}\Big)=0.
\]
Therefore, when $\alpha=\bar{\alpha}$, the matrix $\Mb$ has an eigenvalue at $-1$.

\vspace{-3mm}
\setlength{\parindent}{3ex}
In order to highlight the relation between $\alpha$ and $\Qb(\lambda)$ defined in \eqref{Qb}, we write $\Qb(\lambda)$ as $\Qb(\lambda,\alpha)$. For any fixed $\lambda \in\mathbb{R}$ and fixed {$\alpha>0$}, we denote $\mu(\lambda,\alpha)$ as the smallest eigenvalue of the matrix $ \Qb(\lambda,\alpha)$.
It follows from Rayleigh–Ritz theorem that $\mu(\lambda,\alpha)=\min_{\| \vb \| =1} \vb^\top \Qb(\lambda,\alpha) \vb$.

\vspace{-3mm}
\setlength{\parindent}{3ex}
Next, we show that for any $\alpha>\bar{\alpha}$, there exists $\lambda(\alpha)\le-1$ such that $\mu(\lambda(\alpha),\alpha)=0$.  
Firstly, note that $\mu(\lambda,\alpha)$ is continuous with respect to $(\lambda,\alpha)$ and $\mu(-1,\bar{\alpha})=0$. Secondly, for any fixed $\alpha$, from the definition of $\Qb(\lambda,\alpha)$, we know that there exists $\lambda_0(\alpha)<-1$ such that $\Qb(\lambda,\alpha)$ is positive definite for all $\lambda\le\lambda_0(\alpha)$, i.e., $\mu(\lambda,\alpha)>0$ for all $\lambda\le\lambda_0(\alpha)$. Thirdly, for any fixed $\lambda <-1$, from $\Qb(\lambda,\alpha)=(\Wb-\lambda\Ib_N)^2 \otimes \Ib_m +\alpha(\lambda-1) \tilde{\Hb}$, we know that $\mu(\lambda,\alpha)$ decreases as $\alpha$ increases.
Thus, $\mu(-1,\alpha)\le\mu(-1,\bar{\alpha})=0$ for all $\alpha\ge\bar{\alpha}$. Hence, for any fixed $\alpha>\bar{\alpha}$, $\mu(\lambda,\alpha)$ is continuous with respect to $\lambda$, $\mu(\lambda,\alpha)>0$ for all $\lambda\le\lambda_0(\alpha)$, and $\mu(-1,\alpha)\le0$. Thus, there exists $\lambda(\alpha)\in[\lambda_0(\alpha),-1]$ such that $\mu(\lambda(\alpha),\alpha)=0$.

\vspace{-3mm}
\setlength{\parindent}{3ex}
Finally, we note that $\det \big(\Qb(\lambda,\alpha)\big)=0$ if $\mu(\lambda,\alpha)=0$.
Therefore it follows from \eqref{det-M} that $\lambda(\alpha)\in[\lambda_0(\alpha),-1]$, is an eigenvalue of the matrix $\Mb$, which is less than or equal to $-1$. Hence, the result follows.

\vspace{-3mm}
\section{Proof of Lemma~\ref{lemma-semisimple-directed}} \label{sec-proof-lemma2}
\vspace{-3mm}
\setlength{\parindent}{3ex}
We first show that $1$ is an eigenvalue of the matrix $\Mb$ and its geometric multiplicity is $m$. 
Note that $1$ is an eigenvalue of the matrix $\Mb$, if and only if there is a nonzero right eigenvector $[\vb^\top_1 \; \vb^\top_2]^\top$, such that
\vspace{-3mm}
\begin{equation*}
\begin{bmatrix}
\Pb\otimes \Ib_m & -\alpha \Ib_{Nm} \\
-\tilde{\Hb} \big((\Ib_N-\Pb) \otimes \Ib_m\big) & \Qb\otimes \Ib_m-\alpha \tilde{\Hb}
\end{bmatrix} \begin{bmatrix}
\vb_1 \\
\vb_2
\end{bmatrix}=\begin{bmatrix}
\vb_1 \\
\vb_2
\end{bmatrix}.
\end{equation*}
These equations are equivalent to
\vspace{-3mm}
\begin{align*}
&-\alpha \vb_2 =\big((\Ib_N-\Pb) \otimes \Ib_m\big) \vb_1, \notag \\
& -\tilde{\Hb}\big((\Ib_N-\Pb) \otimes \Ib_m\big) \vb_1=\big((\Ib_N-\Qb) \otimes \Ib_m+\alpha \tilde{\Hb}\big) \vb_2. \notag 
\end{align*}
With a little bit algebra, the above equations are also equivalent to
\vspace{-3mm}
\begin{subequations}\label{eig-one-structure}
\begin{align}
 -\alpha \vb_2 =&\big((\Ib_N-\Pb) \otimes \Ib_m\big) \vb_1, \label{eig-one-structure-eq1} \\
 {\bf 0}_{Nm}=&\big((\Ib_N-\Qb) \otimes \Ib_m\big) \vb_2. \label{eig-one-structure-eq2}
\end{align}
\end{subequations}
From Assumptions~\ref{directed-SC-ass} and~\ref{mixing-weight-PQ-ass}, we know that the matrix $\Pb$ has a unique eigenvalue at $1$ with the right eigenvector ${\bf 1}_N$ and the left eigenvector $\bm{\omega}^\top>0$ (all elements are positive), i.e.,
\vspace{-2mm}
\begin{equation}\label{Peig}
\Pb {\bf 1}_N={\bf 1}_N, \; \bm{\omega}^\top \Pb =\bm{\omega}^\top, \; \bm{\omega}^\top {\bf 1}_N=1.
\vspace{-2mm}
\end{equation}
Therefore, the matrix $\Ib_N-\Pb$ has a unique eigenvalue at zero with the right eigenvector ${\bf 1}_N$ and the left eigenvector $\bm{\omega}^\top>0$. 

\vspace{-3mm}
\setlength{\parindent}{3ex}
Similarly, the matrix $\Qb$ has a unique eigenvalue at $1$ with the right eigenvector $\bm{\mu}>0$ and the left eigenvector ${\bf 1}^\top_N$, i.e.,
\vspace{-2mm}
\begin{equation}\label{Qeig}
\Qb \bm{\mu}=\bm{\mu}, \; {\bf 1}^\top_N \Qb ={\bf 1}^\top_N, \; {\bf 1}^\top_N \bm{\mu}=1.
\vspace{-2mm}
\end{equation}
Thus, the matrix $\Ib_N-\Qb$ has a unique eigenvalue at zero with the right eigenvector $\bm{\mu}>0$ and the left eigenvector ${\bf 1}^\top_N$. 
Therefore, \eqref{eig-one-structure-eq2} is satisfied if $\vb_2={\bf 0}_{Nm}$ or $\vb_2=\bm{\mu} \otimes \cb_1$, where $\cb_1 \in \R^m$ is nonzero. 
Next, we show that $\vb_2=\bm{\mu} \otimes \cb_1$ does not hold by contradiction. 
Suppose that $\vb_2=\bm{\mu} \otimes \cb_1$. 
Then pre-multiplying \eqref{eig-one-structure-eq1} by  $(\bm{\omega} \otimes \cb_1)^\top$ gives $0=-\alpha \|\cb_1\|^2 \bm{\omega}^\top \bm{\mu}$, which does not hold since $\alpha>0$, $\cb_1$ is nonzero, $\bm{\omega}^\top>0$ and $\bm{\mu}>0$. 
Therefore, $\vb_2={\bf 0}_{Nm}$. 
This together with \eqref{eig-one-structure-eq1}, the above properties of the matrix $\Ib_N-\Pb$ and the fact that $[\vb^\top_1 \; \vb^\top_2]^\top$ is nonzero implies that $\vb_1={\bf 1}_N \otimes \cb$ for some nonzero $\cb \in \R^m$.

\vspace{-3mm}
\setlength{\parindent}{3ex}
Hence, we conclude that $1$ is an eigenvalue of the matrix $\Mb$ with the geometric multiplicity being $m$ and the corresponding right eigenvectors are of the form $[\vb^\top_1 \; \vb^\top_2]^\top$ with $\vb_1={\bf 1}_N \otimes \cb$ for some nonzero $\cb \in \R^m$ and $\vb_2={\bf 0}_{Nm}$.

\vspace{-3mm}
\setlength{\parindent}{3ex}
Next,  we show that the algebraic multiplicity associated with the eigenvalue at $1$ is also $m$ by contradiction.
Suppose that the algebraic multiplicity of the eigenvalue at $1$ is strictly greater than $m$.
Then there exists a nonzero right generalized eigenvector $[\ub^\top_1 \; \ub^\top_2]^\top$, such that
\vspace{-3mm}
\[
(\Mb-\Ib_{Nm}) \begin{bmatrix}
\ub_1 \\
\ub_2
\end{bmatrix}=\begin{bmatrix}
\vb_1 \\
\vb_2
\end{bmatrix},
\vspace{-3mm}
\]
where $[\vb^\top_1 \; \vb^\top_2]^\top$ is the right eigenvector corresponding to the eigenvalue at $1$, i.e., $\vb_1={\bf 1}_N \otimes \cb$ for some nonzero $\cb \in \R^m$ and $\vb_2={\bf 0}_{Nm}$.
This is equivalent to
\vspace{-3mm}
\begin{subequations}
\begin{align*}
&\big((\Pb-\Ib_N) \otimes \Ib_m\big) \ub_1-\alpha \ub_2 ={\bf 1}_N \otimes \cb,\\
& -\tilde{\Hb}\big((\Ib_N-\Pb) \otimes \Ib_m\big) \ub_1+\big((\Qb-\Ib_N) \otimes \Ib_m-\alpha \tilde{\Hb}\big) \ub_2=0.
\end{align*}
\end{subequations}
From the above equations, we obtain
\vspace{-3mm}
\[
\big((\Qb-\Ib_N) \otimes \Ib_m\big) \ub_2=-\tilde{\Hb}({\bf 1}_N \otimes \cb).
\]
We then have
\vspace{-3mm}
\begin{align*}
&({\bf 1}_N \otimes \cb)^\top\tilde{\Hb}({\bf 1}_N \otimes \cb)\\
&=-({\bf 1}_N \otimes \cb)^\top\big((\Qb-\Ib_N) \otimes \Ib_m\big) \ub_2=0,
\vspace{-3mm}
\end{align*}
where the second equality follows from the fact that ${\bf 1}^\top_N \Qb={\bf 1}^\top_N$. 
It then follows from Lemma~\ref{lemmaH} that $\cb={\bf 0}_m$, which contradicts with the fact that $\cb$ is a nonzero vector.
Therefore, the matrix $\Mb$ always has $m$ semisimple eigenvalues at $1$.

\vspace{-3mm}
\section{Proof of Theorem~\ref{directed-thm2}}\label{sec-proof-directed}
\vspace{-3mm}
\setlength{\parindent}{3ex}
In view of Proposition~\ref{directed-thm1}, in order to prove the theorem, it is equivalent to show that all the
eigenvalues of the matrix $\Mb$, except $m$ semisimple eigenvalues at $1$, are strictly within the unit circle if $\alpha$ is sufficiently small.

\vspace{-3mm}
\setlength{\parindent}{3ex}
In order to show this, we use the eigenvalue perturbation theory.
Note that $\mathbf{M}=\mathbf{M}_0+\alpha {\bm \Delta}$, where
\vspace{-2mm}
\begin{align*}
\mathbf{M}_0=&\begin{bmatrix}
\Pb\otimes \Ib_m & {\bf 0}_{Nm \times Nm} \\
-\tilde{\Hb} \big((\Ib_N-\Pb) \otimes \Ib_m\big) & \Qb\otimes \Ib_m
\end{bmatrix},\\
\bm{\Delta} =& \begin{bmatrix}
 {\bf 0}_{Nm \times Nm}  & -\Ib_{Nm} \\
 {\bf 0}_{Nm \times Nm} & -\tilde{\Hb}
\end{bmatrix}.
\end{align*}
Thus the matrix $\Mb$ can be viewed as the matrix $\Mb_0$ perturbed by $\alpha \bm{\Delta}$.

\vspace{-3mm}
\setlength{\parindent}{3ex}
Note that the from Assumptions~\ref{directed-SC-ass} and~\ref{mixing-weight-PQ-ass}, we know that the matrices $\Pb$ and $\Qb$ has a unique eigenvalue at $1$, and all other eigenvalues are strictly within the unit circle. 
Given the block diagonal structure of the matrix $\Mb_0$, it is easy to see that the matrix $\Mb_0$ has $2m$ semisimple eigenvalues at $1$, and all other eigenvalues are strictly within the unit circle.
Denote $\lambda_i$, $i=1,\ldots,2Nm$ as the eigenvalues of the matrix $\Mb_0$.
Without loss of generality, assume that $\lambda_j=1$ for $j=1,\ldots, 2m$.
It is easy to verify that the right eigenvectors associated with the  eigenvalues at $1$ of the matrix $\mathbf{M}_0$ are
\vspace{-3mm}
\begin{equation*}
\begin{aligned}
\Ub \triangleq &\begin{bmatrix}
{\bf 0}_N \otimes \Ib_m &&  {\bf 1}_N \otimes \Ib_m \\
\bm{\mu} \otimes \Ib_m & & (\bm{\mu}  \otimes \Ib_m) ({\bf 1}^\top_N \otimes \Ib_m) \tilde{\Hb} ({\bf 1}_N \otimes \Ib_m)
\end{bmatrix},
\end{aligned}
\vspace{-2mm}
\end{equation*}
and the corresponding left eigenvectors are
\vspace{-3mm}
\[
\bm{\Omega}^\top \triangleq
\begin{bmatrix}
-({\bf 1}_N^\top \otimes \Ib_m) \tilde{\Hb} & &  {\bf 1}_N^\top \otimes \Ib_m \\
\bm{\omega}^\top \otimes \Ib_m & & {\bf 0}_N^\top \otimes \Ib_m
\end{bmatrix},
\vspace{-2mm}
\]
where $\bm{\omega}^\top>0$ (all elements are positive) and $\bm{\mu}>0$ are given by \eqref{Peig} and \eqref{Qeig}, respectively.

\vspace{-3mm}
\setlength{\parindent}{3ex}
Denote $\lambda_i(\alpha)$, $i=1,\ldots,2Nm$ as the eigenvalues of the matrix $\Mb=\Mb_0+\alpha \bm{\Delta}$ corresponding relatively to $\lambda_i$. 
Note that $\bm{\Omega}^\top \Ub=\Ib_{2m}$. 
It then follows from the eigenvalue perturbation theory \citep[Lemma~7]{cai-kai-surplus-fixed} that
the derivatives $\frac{d}{d \alpha}\lambda_j(\alpha)|_{\alpha=0}$ for $j=1,\ldots, 2m$ exist and are the eigenvalues of the following matrix:
\vspace{-2mm}
\[
\bm{\Omega}^\top \bm{\Delta} {\Ub} = \begin{bmatrix}
{\bf 0}_{m \times m} & {\bf 0}_{m \times m} \\
-\bm{\omega}^\top \bm{\mu} \Ib_m & -\bm{\omega}^\top \bm{\mu}({\bf 1}^\top_N \otimes \Ib_m) \tilde{\Hb} ({\bf 1}_N \otimes \Ib_m)
\end{bmatrix}.
\]

\vspace{-3mm}
\setlength{\parindent}{3ex}
Moreover, due to the block triangular structure, the above matrix has $m$ eigenvalues at $0$, and the other $m$ eigenvalues are the eigenvalues of the matrix  $-\bm{\omega}^\top \bm{\mu}({\bf 1}^\top_N \otimes \Ib_m) \tilde{\Hb} ({\bf 1}_N \otimes \Ib_m)$.
Given the structures of $\Hb$ and $\tilde{\Hb}$,
we obtain that
\vspace{-2mm}
\begin{equation}\label{sum-pd}
({\bf 1}^\top_N \otimes \Ib_m) \tilde{\Hb} ({\bf 1}_N \otimes \Ib_m)=\Hb^\top \Hb,
\vspace{-2mm}
\end{equation}
which is positive definite since $\rank(\Hb)=m$ from Assumption~\ref{unique-LS-sol-ass}.
Also note that $-\bm{\omega}^\top \bm{\mu}<0$ since $\bm{\omega}^\top>0$ and $\bm{\mu}>0$.
Therefore, all eigenvalues of the matrix $-\bm{\omega}^\top \bm{\mu}({\bf 1}^\top_N \otimes \Ib_m) \tilde{\Hb} ({\bf 1}_N \otimes \Ib_m)$ are negative.

\vspace{-3mm}
\setlength{\parindent}{3ex}
Thus, $\frac{d}{d \alpha}\lambda_j(\alpha)|_{\alpha=0}=0$ for $j=1,\ldots, m$ and $\frac{d}{d \alpha}\lambda_j(\alpha)|_{\alpha=0}<0$ for $j=m+1, \ldots, 2m$.
Since the eigenvalues of the matrix $\mathbf{M}$ are continuous of $\alpha$,
as $\alpha$ slightly increases from zero, the eigenvalues $\lambda_1(\alpha), \ldots, \lambda_m(\alpha)$ stay at $1$, while  $\lambda_{m+1}(\alpha), \ldots, \lambda_{2m}(\alpha)$ move to the left along the real axis.
Let $\delta_1$ be the upper bound of $\alpha$ such that when $\alpha<\delta_1$, $|\lambda_{i}(\alpha)|<1$ for $i=m+1, \ldots, 2m$.
Since the eigenvalues of the matrix $\Mb=\Mb_0+\alpha \bm{\Delta}$ continuously depend on $\alpha$,
there exists an upper bound $\delta_2$ such that when $\alpha<\delta_2$, $|\lambda_{i}(\alpha)|<1$ for $i=2m+1, \ldots, 2Nm$.
Hence, for sufficiently small $\alpha<\alpha_1=\min\{\delta_1,\delta_2\}$, the matrix $\Mb$ has $m$ eigenvalues at $1$ while all other eigenvalues are strictly within the unit circle.


\begin{thebibliography}{38}
	
	
	\bibitem[{Anderson et~al.(2016)Anderson, Mou, Morse, and
		Helmke}]{anderson2015decentralized}
	Anderson, B., Mou, S., Morse, A.~S., Helmke, U., 2016. Decentralized gradient
	algorithm for solution of a linear equation. Numerical Algebra, Control and
	Optimization 6~(3), 319--328.
	
	\bibitem[{Arrow et~al.(1958)Arrow, Huwicz, and Uzawa}]{Arrow58}
	Arrow, K., Huwicz, L., Uzawa, H., 1958. Studies in Linear and Non-linear
	Programming. Stanford University Press.
	
	\bibitem[{Cai and Ishii(2012)}]{cai-kai-surplus-fixed}
	Cai, K., Ishii, H., 2012. Average consensus on general strongly connected
	digraphs. Automatica 48~(11), 2750 -- 2761.
	
	\bibitem[{Charalambous et~al.(2015)Charalambous, Yuan, Yang, Pan, Hadjicostis,
		and Johansson}]{Themis-TCNS}
	Charalambous, T., Yuan, Y., Yang, T., Pan, W., Hadjicostis, C.~N., Johansson,
	M., 2015. Distributed finite-time average consensus in digraphs in the
	presence of time-delays. {IEEE} Transactions on Control of Network Systems
	2~(4), 370--381.
	
	\bibitem[{Cort\'{e}s(2008)}]{Minimum-consensus}
	Cort\'{e}s, J., 2008. Distributed algorithms for reaching consensus on general
	functions. Automatica 44~(3), 726--737.
	
	\bibitem[{Du et~al.(2018)Du, Yao, Wu, Li, Liu, and Yang}]{Yang-PESGM2018}
	Du, W., Yao, L., Wu, D., Li, X., Liu, G., Yang, T., 2018. Accelerated
	distributed energy management for microgrids. In: Proc. of the IEEE Power and
	Energy Society General Meeting.
	
	\bibitem[{Gharesifard and Cort\'{e}s(2014)}]{Cortes_TAC_CT}
	Gharesifard, B., Cort\'{e}s, J., 2014. Distributed continuous-time convex
	optimization on weight-balanced digraphs. {IEEE} Transactions on Automatic
	Control 59~(3), 781--786.
	
	\bibitem[{Horn and Johnson(2001)}]{Horn2001}
	Horn, R.~A., Johnson, C.~R., 2001. Matrix Analysis, 2nd Edition. Cambridge
	University Press.
	
	\bibitem[{Jakoveti{\'c}(2019)}]{Dusan-Unified}
	Jakoveti{\'c}, D., 2019. A unification and generalization of exact distributed
	first-order methods. IEEE Transactions on Signal and Information Processing
	over Networks 5~(1), 31--46.
	
	\bibitem[{Jury(1991)}]{jury}
	Jury, E., 1991. A note on the modified stability table for linear discrete time
	systems. IEEE Trans. Circ. \& Syst. 38~(2), 221--223.
	
	\bibitem[{Liu et~al.(2017)Liu, Morse, Nedi\'{c}, and Ba\c{s}ar}]{JiLiu-Auto}
	Liu, J., Morse, A.~S., Nedi\'{c}, A., Ba\c{s}ar, T., 2017. Exponential
	convergence of a distributed algorithm for solving linear algebraic
	equations. Automatica 83, 37--46.
	
	\bibitem[{Liu et~al.(2018)Liu, Mou, and Morse}]{JiLiu-TAC}
	Liu, J., Mou, S., Morse, A.~S., 2018. Asynchronous distributed algorithms for
	solving linear algebraic equations. {IEEE} Transactions on Automatic Control
	63~(2), 372--385.
	
	\bibitem[{Liu et~al.(2019)Liu, Lageman, Anderson, and Shi}]{Guodong-arXiv}
	Liu, Y., Lageman, C., Anderson, B. D.~O., Shi, G., 2019. An
	{A}rrow-{H}urwicz-{U}zawa type flow as least squares solver for network
	linear equations. Automatica 100, 187--193.
	
	\bibitem[{Lu and Tang(2018)}]{Jie-TCNS2018}
	Lu, J., Tang, C.~Y., 2018. A distributed algorithm for solving positive
	definite linear equations over networks with membership dynamics. {IEEE}
	Transactions on Control of Network Systems 5~(1), 215--227.
	
	\bibitem[{Mou et~al.(2015)Mou, Liu, and Morse}]{Mou-TAC15}
	Mou, S., Liu, J., Morse, A.~S., 2015. A distributed algorithm for solving a
	linear algebraic equation. {IEEE} Transactions on Automatic Control 60~(11),
	2863--2878.
	
	\bibitem[{Mou et~al.(2016)Mou, Morse, Lin, Wang, and Fullmer}]{Mou_SCL16}
	Mou, S., Morse, A.~S., Lin, Z., Wang, L., Fullmer, D., 2016. A distributed
	algorithm for efficiently solving linear equations and its applications.
	Systems \& Control Letters 91, 21--27.
	
	\bibitem[{Nedi\'{c} et~al.(2017)Nedi\'{c}, Olshevsky, and Shi}]{Nedic-arXiv1}
	Nedi\'{c}, A., Olshevsky, A., Shi, W., 2017. Achieving geometric convergence
	for distributed optimization over time-varying graphs. SIAM Journal on
	Optimization 27~(4), 2597--2633.
	
	\bibitem[{Nedi\'{c} and Ozdaglar(2009)}]{Nedic09}
	Nedi\'{c}, A., Ozdaglar, A., 2009. Distributed subgradient methods for
	multi-agent optimization. {IEEE} Transactions on Automatic Control 54~(1),
	48--61.
	
	\bibitem[{Nedi\'{c} et~al.(2010)Nedi\'{c}, Ozdaglar, and Parrilo}]{Nedic10}
	Nedi\'{c}, A., Ozdaglar, A., Parrilo, P.~A., 2010. Constrained consensus and
	optimization in multi-agent networks. {IEEE} Transactions on Automatic
	Control 55~(4), 922--938.
	
	\bibitem[{Olfati-Saber et~al.(2007)Olfati-Saber, Fax, and Murray}]{Olfati07}
	Olfati-Saber, R., Fax, J.~A., Murray, R.~M., 2007. Consensus and cooperation in
	networked multi-agent systems. Proceedings of the {IEEE} 95~(1), 215--233.
	
	\bibitem[{Pu et~al.(2018)Pu, Shi, Xu, and Nedi\'{c}}]{pu2018push}
	Pu, S., Shi, W., Xu, J., Nedi\'{c}, A., 2018. A push-pull gradient method for
	distributed optimization in networks. In: 57th IEEE Conference on Decision
	and Control (CDC). pp. 3385--3390.
	
	\bibitem[{Qu and Li(2018)}]{Lina-TCNS}
	Qu, G., Li, N., 2018. Harnessing smoothness to accelerate distributed
	optimization. {IEEE} Transactions on Control of Network Systems 5~(3),
	1245--1260.
	
	\bibitem[{Shi et~al.(2017)Shi, Anderson, and Helmke}]{Guodong-TAC17}
	Shi, G., Anderson, B. D.~O., Helmke, U., 2017. Network flows that solve linear
	equations. {IEEE} Transactions on Automatic Control 62~(6), 2659--2674.
	
	\bibitem[{Shi et~al.(2015)Shi, Ling, Wu, and Yin}]{Yin-EXTRA}
	Shi, W., Ling, Q., Wu, G., Yin, W., 2015. {EXTRA}: An exact first-order
	algorithm for decentralized consensus optimization. SIAM Journal on
	Optimization 25~(2), 944--966.
	
	\bibitem[{Silvester(2000)}]{Silvester-det}
	Silvester, J.~R., 2000. Determinants of block matrices. The Mathematical
	Gazette 84~(501), 460--467.
	
	\bibitem[{Sundaram and Hadjicostis(2007)}]{Sunda07}
	Sundaram, S., Hadjicostis, C.~N., 2007. Finite-time distributed consensus in
	graphs with time-invariant topologies. In: Proc. American Control Conference.
	pp. 711--716.
	
	\bibitem[{Wang and Elia(2010)}]{Elia_Allerton10}
	Wang, J., Elia, N., 2010. Control approach to distributed optimization. In:
	Proc. 48th Annual Allerton Conference on Communication, Control, and
	Computing (Allerton). pp. 557--561.
	
	\bibitem[{Wang and Elia(2012)}]{Elia_ACC2012}
	Wang, J., Elia, N., 2012. Distributed least square with intermittent
	communications. In: Proc. American Control Conference. pp. 6479--6484.
	
	\bibitem[{Wang et~al.(2019{\natexlab{a}})Wang, Ren, and Duan}]{PengWang-TCNS}
	Wang, P., Ren, W., Duan, Z., 2019{\natexlab{a}}. Distributed algorithm to solve
	a system of linear equations with unique or multiple solutions from arbitrary
	initializations. {IEEE} Transactions on Control of Network Systems 6~(1),
	82--93.
	
	\bibitem[{Wang et~al.(2019{\natexlab{b}})Wang, Zhou, Mou, and
		Corless}]{Mou-TAC}
	Wang, X., Zhou, J., Mou, S., Corless, M.~J., 2019{\natexlab{b}}. A distributed
	algorithm for least squares solutions. {IEEE} Transactions on Automatic
	Control, to appear.
	
	\bibitem[{Xiao and Boyd(2004)}]{Boyd_SCL04}
	Xiao, L., Boyd, S., 2004. Fast linear iterations for distributed averaging.
	Systems \& Control Letters 53~(1), 65--78.
	
	\bibitem[{Xin and Khan(2018)}]{Usman-CSL2018}
	Xin, R., Khan, U.~A., 2018. A linear algorithm for optimization over directed
	graphs with geometric convergence. Control Systems Letters 2~(3), 315--320.
	
	\bibitem[{Xu et~al.(2015)Xu, Zhu, Soh, and Xie}]{Xu_CDC15}
	Xu, J., Zhu, S., Soh, Y.~C., Xie, L., 2015. Augmented distributed gradient
	methods for multi-agent optimization under uncoordinated constant stepsizes.
	In: Proc. of the 54th IEEE Conference on Decision and Control. pp.
	2055--2060.
	
	\bibitem[{Xu et~al.(2018)Xu, Zhu, Sohy, and Xie}]{Lihua-Proximal-TAC}
	Xu, J., Zhu, S., Sohy, Y.~C., Xie, L., 2018. A {B}regman splitting scheme for
	distributed optimization over networks. {IEEE} Transactions on Automatic
	Control 63~(11), 3809--3824.
	
	\bibitem[{Yang et~al.(2016)Yang, Wu, Sun, and Lian}]{YangT16}
	Yang, T., Wu, D., Sun, Y., Lian, J., 2016. Minimum-time consensus based
	approach for power system applications. IEEE Transactions on Industrial
	Electronics 63~(2), 1318--1328.
	
	\bibitem[{Yao et~al.(2018)Yao, Yuan, Sundaram, and Yang}]{LishaoYao-ICCA2018}
	Yao, L., Yuan, Y., Sundaram, S., Yang, T., 2018. Distributed finite-time
	optimization. In: Proc. 14th IEEE International Conference on Control \&
	Automation. pp. 147--154.
	
	\bibitem[{Yuan et~al.(2013)Yuan, Stan, Shi, Barahona, and Goncalves}]{Ye-AUTO}
	Yuan, Y., Stan, G.-B., Shi, L., Barahona, M., Goncalves, J., 2013.
	Decentralised minimum-time consensus. Automatica 49~(5), 1227--1235.
	
	\bibitem[{Zeng et~al.(2019)Zeng, Liang, Hong, and Chen}]{Hong-TAC18}
	Zeng, X., Liang, S., Hong, Y., Chen, J., 2019. Distributed computation of
	linear matrix equations: An optimization perspective. {IEEE} Transactions on
	Automatic Control 64~(5), 1858--1873.
\end{thebibliography}
\end{document}